\documentclass[12pt]{article}
\usepackage{amssymb,latexsym,amsmath}
\def\C{\mathbb C}
\def\R{\mathbb R}
\def\N{\mathbb N}
\def\Z{\mathbb Z}
\def\Q{\mathbb Q}
\numberwithin{equation}{section}

\newtheorem{thm}{Theorem}[section]
\newtheorem{lem}{Lemma}[section]

\newtheorem{prop}{Proposition}[section]

\begin{document}
\sffamily
\title{Linear differential polynomials in\\ zero-free meromorphic functions}
\author{J.K. Langley}
\maketitle
\centerline{\textit{In fondest memory of G\"unter Frank and Milne Anderson}}

\begin{abstract}
The paper determines all meromorphic functions $f$ in $\C$ such that
 $f $ and $F$  have finitely many zeros,
 where $F = f^{(k)} + a_{k-1} f^{(k-1)} + \ldots + a_0f$ with $k \geq 3$ and
the  $a_j$ rational functions. 
\end{abstract}

\section{Introduction}

Let the function $f$ be meromorphic 
in an annulus $\Omega (r_1) = \{ z \in \C : r_1 < |z| < \infty \} $, with $r_1$ positive
(not necessarily the same at each occurrence in this paper). 
Let $k \geq 2$ and  let $a_0, \ldots , a_{k-1}$ be functions which are rational at infinity, that is,
analytic on some  $\Omega(r_1)$ with at most a pole at $\infty$. 
Write $D = d/dz$ and 
\begin{equation}
 \label{1}
F = L[f], \quad L = D^k + a_{k-1}D^{k-1}  + \ldots +  a_0 ,
\end{equation}
in which  $L[y]$ denotes the operator $L$ acting on the function $y$. 
The central objective of
this paper is the classification
of all those $f$ for which 
$f$ and $F$ have no zeros in $\Omega(r_1)$. 
By a standard change of variables $f = e^Pg, F = e^P G$, with $P$ a polynomial,
it may be assumed that $a_{k-1}( \infty ) = 0$. 

This problem, part of which appeared as 1.42 in the
collection \cite{branhay}, has a long history going
back to Hayman's conjecture in \cite{Hay1}, proved 
in \cite{Fra1,La5}, that if $k \geq 2$ then the only 
meromorphic
functions $f$ in the plane for which $f$ and $f^{(k)}$ have no zeros are those of form 
$f(z) = e^{az+b}$ or $f(z) = (az+b)^{-n}$ with $a, b \in \C$ and $n \in \N = \{ 1, 2, \ldots \}$: more generally, if 
$f$ and $f^{(k)}$ have finitely many zeros then $f = Re^P$, 
with $R$ a rational function and $P$ 
a polynomial \cite{FHP,La5}, so that $f'/f$ is rational. 
The  problem for $k=2$ and coefficients which are 
rational at infinity was fully solved in \cite{La5,La9}. 

\begin{thm}[\cite{La5,La9}]
 \label{thmA}
Let the function $f$ be meromorphic in $S \leq |z| < \infty$ for some $S > 0$
and let the functions $a _1$ and $a_0$ be analytic there and rational at infinity. Assume that
$a_1(\infty) = 0$ and that $f$ and $F=f''+a_1f' +a_0f$ have no zeros in $S \leq |z| < \infty$. 
\\
(a) If 
$$
{\rm deg}_\infty (a_0) = \lim_{z \to \infty} \frac{ \log |a_0(z)|}{\log |z|} 
$$
is even then at least one of the following holds.\\
(i) The function $f'/f$ is rational at infinity.\\
(ii) The function $f$ satisfies 
\begin{equation}
 \label{rdm1}
\frac{f'}{f} = - \frac{a_1}2 +    \frac{g'}{2 g} + \frac{A}{g}  ,  \quad 
g^2 = \frac{f}{F} , \quad 
g' =  \left( 2 \frac{f_1'}{f_1} + a_1 \right) g + B, 
\end{equation}
where $A, B \in \C$ and $g$ is analytic in $|z| \geq S$, while  
$f_1$ is a solution of the homogeneous equation 
\begin{equation}
 \label{rdm3}
w''+a_1 w' + a_0w = 0 
\end{equation}
which admits unrestricted analytic continuation without zeros in $|z| \geq S $.\\
(iii) There exist solutions $f_1$, $f_2$ of (\ref{rdm3}), such that
\begin{equation}
 \label{rdm4}
f = A f_2 \left( 1 + B \left( \frac{f_2}{f_1} \right)^{1/N} \right)^{-N} , 
\quad A, B \in \C, \quad N \in \N.
\end{equation}
Here both $f_1$ and $f_2$ admit unrestricted analytic continuation without zeros in $|z| > R_1$ for some $R_1 > 0$,
and $(f_2/f_1)^{1/N} $ is analytic in $|z| > R_1$.\\
(iv) There exist solutions $f_1$, $f_2$ of (\ref{rdm3}), each admitting unrestricted analytic continuation 
without zeros in $|z| > R_1$ for some $R_1 > 0$, a function $M$ which is rational at infinity, 
and non-constant polynomials $Q$, $Q_1$ such that
\begin{equation*}
 \label{rdm5}
\frac{f'}{f} = \frac{f_2'}{f_2} + \frac{Q(M) M'}{e^M + 1} ,
\quad \hbox{where} \quad 
Q(M)M' = \frac{f_1'}{f_1} - \frac{f_2'}{f_2}  \quad \hbox{or} \quad Q_1(M)e^{-M} = \frac{f_1}{f_2} .
\end{equation*}
(b) If  ${\rm deg}_\infty (a_0)$ is  odd then $f$ may be determined by applying part (a) to 
$$\phi(z) = f(z^2), \quad  \Phi(z) = 4z^2 F(z^2) = \phi''(z) + (2za_1(z^2)-1/z) \phi'(z) + 4z^2 a_0(z^2) \phi(z) .
$$
\end{thm}

A refinement of this theorem for  meromorphic functions in the plane may be found in \cite[Theorem 1.3]{La2013}. 
For $k \geq 3$ and $f, F$ zero-free in the whole plane, the case of 
constant coefficients was solved in full by Steinmetz  in \cite{Stei1}, while
polynomial coefficients were treated in \cite{FH}
for entire $f$,
and for meromorphic $f$ by Br\"uggemann in \cite{Bru}.

\begin{thm}
[\cite{Bru,FH}] 
 \label{thmBru}
Let the function $f$ be meromorphic in the plane, such that $f$ and $F = L[f]$ have no zeros, where $k \geq 3$ and
$a_0, \ldots, a_{k-2}$ are polynomials, not all constant, with 
$a_{k-1} \equiv 0$. Then $f = (H')^{-(k-1)/2} e^H $ or $f = (H')^{-(k-1)/2} H^{-m} $ 
for some $m \in \N$, where $H''/H'$ is a polynomial. 
\end{thm}

The following theorem, which settles all cases, will be proved. 

\begin{thm}
 \label{thm1}
Let $k \geq 3$ and let the function $f$ be meromorphic in some annulus $\Omega(r_1)$, with $f'/f$ not rational at infinity.
Assume that $f$ and $F=L[f]$ have no zeros in $\Omega(r_1)$,  where $L$ is as in (\ref{1})
with the $a_j$ analytic in $\Omega(r_1)$ and rational at infinity,  and with
$a_{k-1}(\infty)
= 0$. Then $f$ satisfies at least one of the following. \\
(i) The logarithmic derivative $f'/f$ has a representation 
\begin{equation}
 \frac{f'}{f} = - \frac{a_{k-1}}k - \left( \frac{k-1}2 \right)  \frac{H''}{H'} + H' 
\quad \hbox{or} \quad 
\frac{f'}{f} = - \frac{a_{k-1}}k - \left( \frac{k-1}2 \right)  \frac{H''}{H'} - m \frac{H'}{H} ,
\label{Brurep}
\end{equation}
where $m \in \N$ and $H_0 = H''/H'$ is rational at infinity, with $H_0(\infty) \neq 0$, while the equation $L[y]=0$ has linearly independent local solutions $y_j$ satisfying
\begin{equation}
 \label{Brurep2}
\frac{y_j'}{y_j} = - \frac{a_{k-1}}k - \left( \frac{k-1}2 \right)  \frac{H''}{H'} + (j-1)  \frac{H'}{H} , \quad j=1, \ldots, k,
\end{equation}
and $f$ is given locally by either $f = cy_1 \exp( y_2/y_1) $ or $f = cy_1^{m+1} y_2^{-m}$, where $c \in \C \setminus \{ 0 \}$.
\\
(ii) There exist a polynomial $Q$ and functions $\nu_1$, $\nu_0$, both rational at infinity, such that $f'/f$ has a representation
\begin{equation}
 \label{caseIIrep}
\frac{f'}{f} = \frac{Q(T) T'}{1-e^{-T} }  + \frac{y_1'}{y_1} , \quad T = \log
\left( \frac{v}{u} \right) , 
\end{equation}
where $y_1$ is a solution of $L[y]=0$, while $v$ and $u$ are linearly independent solutions of 
\begin{equation}
 \label{b1}
y'' + \nu_1 y' + \nu_0 y =0
\end{equation}
which continue without zeros in some annulus $\Omega(r_2)$. 
Here $Q(T)$  is rational at infinity, and $u$, $v$,  $y_1'/y_1$ and $a_0, \ldots, a_{k-2}$
all have representations in terms of
$Q(T)$, $T$, $a_{k-1}$ and their derivatives. 
Moreover, if $T'$ is  not rational at infinity then $k$ is even and 
$z^{-1/2} T'(z)$ is rational at infinity.

In both cases (i) and (ii) 
there exist $r_3 > 0$ and  functions $\widetilde a_1$, $\widetilde a_0$, each rational at infinity,
such that $f'' + \widetilde a_1 f' + \widetilde a_0 f$ has no zeros in  $\Omega(r_3)$. 
\end{thm}

The conclusions 
of Theorems \ref{thmA} and \ref{thm1} are closely related, and the last assertion of Theorem~\ref{thm1} makes it clear that
this is no coincidence. 
If $Q$ is a constant $d$ in (\ref{caseIIrep}) then integration shows that $f$ is a constant multiple of $y_1 \left(  v/u -1 \right)^{d} $. Conclusion (\ref{Brurep})  
may be compared with that of Theorem~\ref{thmBru}, and links closely to (\ref{rdm1}) of Theorem \ref{thmA} and
\cite[Theorem 1.3(II)]{La2013}. Examples II and III in Section \ref{exa} demonstrate that in (\ref{caseIIrep}) the multiplicities
of poles of $f$ may be unbounded, in sharp contrast to the situation in Theorem \ref{thmBru}, where any poles of $f$ must all have the same multiplicity $m$.
Example III also shows  that $T'$ need not be rational at infinity in (\ref{caseIIrep}).

Some previous partial results for rational coefficients may be found in \cite{Hec2,Latsuji}. 
Methods from \cite{Bru,Fra1,FH,Stei1} are essential to the proof of Theorem \ref{thm1}; these are supplemented 
by a result (Lemma \ref{lemintegervalued}) on integer-valued analytic functions, facilitating the analytic continuation of several asymptotic representations.
A decisive role is played by a criterion (Lemma \ref{lemtwogsmall}) for certain auxiliary functions to satisfy a second order differential equation,
which simplifies the subsequent analysis considerably.

The author   acknowledges extensive discussions and correspondence on this problem 
with the late
G\"unter Frank; these took place over many years and have contributed 
substantially to the methodology of this paper. Indeed, the 
Wronskian-based method invented by
Frank \cite{Fra1,FHP} underpins much of the 
successful work on these and related problems.
Thanks are also due to the referees  for their valuable comments. 

\section{Examples}\label{exa}

Throughout the paper  $c$ will be used to denote
non-zero constants, not always the same at each occurrence, and
$\C^*$ will denote $\C \setminus \{ 0 \}$. 

\subsection{Example I}

This example goes back to \cite{FH}, and may be compared with  conclusion (i)
of Theorem \ref{thm1} and that of Theorem \ref{thmBru}.  
Let $H$ be such that $\delta = H''/H' \not \equiv 0$ is  a polynomial, and write
$$
g = (H')^{-k} e^H, \quad h = (H')^{-k} H^{-m} , \quad  D = \frac{d\,}{dz} , \quad m \in \N .
$$
Then it is 
easy to check (see the remark following (\ref{Brurep5}) below) that 
$$
 (D + \delta ) \ldots (D + k \delta ) [g ] = e^H, \quad (D + \delta ) \ldots (D + k \delta ) [h] = c H^{-m-k} .
$$
Taking $f$ to be $e^P g$ or $e^P h$ for a suitably chosen polynomial $P$ 
gives polynomial coefficients $a_j$ with $a_{k-1}=0$ such that $f$ and $F=L[f]$ have no zeros. 

\subsection{Example II}\label{ex2}

Let $P $ be a  non-constant polynomial 
which takes positive integer values at all zeros of
$1- e^z $, and write 
\begin{equation}
 \label{deta1}
\frac{f'(z)}{f(z)} = \frac{P(z)}{1- e^z }  , \quad 
\frac{f''(z)}{f(z)} 
= \frac{Q_1(z)e^z + Q_0(z)}{(1-e^z)^2} , \quad 
Q_1 = P - P', \quad Q_0 = P'+ P^2 .
\end{equation}
Then $f$ is meromorphic and zero-free in the plane, with a pole of multiplicity $P(z)$ at a zero $z$
of $1-e^z$. A standard calculation yields
polynomials $R_j$ such that 
$$
\frac{f'''(z)}{f(z)} =  \frac{R_2(z) e^{2z} + R_1(z)e^z + R_0(z)}{(1-e^z)^3} .
$$
If $F = f'''+ b_2 f'' + b_1 f' $, where the $b_j$ are rational functions, then 
\begin{eqnarray*}
  \frac{F(z)}{f(z)} 
&=& \frac{ B_2(z) e^{2z} + B_1(z) e^z + B_0(z)}{(1-e^z)^3} ,\\
B_2 &=& R_2 - b_2 Q_1 + b_1 P , \\
B_1 &=& R_1 + b_2(Q_1-Q_0) - 2b_1 P ,  \\
B_0 &=&  R_0 + b_2 Q_0 + b_1 P  .
\end{eqnarray*}
Thus $F$ may be made zero-free  in some $\Omega(r_1)$ by setting 
\begin{eqnarray}
 \label{exL1}
0 &=&   R_1 + b_2(Q_1-Q_0) - 2b_1 P  =    R_0 + b_2 Q_0 + b_1 P  , \nonumber \\ 
\frac{F(z)}{f(z)} &=& \frac{B_2(z) e^{2z}}{(1-e^z)^3} ,
\end{eqnarray}
these equations being solvable for $b_1$ and $b_2$, since $(Q_1-Q_0)P+ 2Q_0 P  = (Q_1+Q_0) P \not \equiv 0$ by (\ref{deta1}).
Similar  calculations show that it is possible to achieve each of 
\begin{equation}
 \label{exL2}
\frac{F(z)}{f(z)} = \frac{B_1(z)e^z }{(1-e^z)^3} \quad ;  \quad  \frac{F(z)}{f(z)} = \frac{B_0(z) }{(1-e^z)^3} .
\end{equation}
Finally, should it be the case that 
$b_2(\infty) \neq 0$, there exist a polynomial $Q_2$ and rational functions $a_j$, with $a_2(\infty) = 0$, such that 
writing $h = e^{Q_2} f $ gives 
$$
\frac{F(z)}{f(z)} = 
\frac{ f'''+ b_2 f'' + b_1 f' }{f} = \frac{ h'''+ a_2 h'' + a_1 h' + a_0 h}{h}  .
$$

\subsection{Example III}\label{coshex}

This is adapted from 
\cite{La9}. Let $Y(z) = z^{m/2}$, where $m \in \N$,
and set $h= \cosh Y$. Then $h $ is entire with only simple zeros. Let $P_1$ be an even polynomial which takes negative integer values at all odd
integer multiples of $\pi i /2$, and set $P = P_1(Y)$.
Then $P$ is a polynomial and setting 
$$
\frac{f'}f = P \cdot \frac{h'}{h} = \frac{PY' \sinh Y}{\cosh Y} = \frac{-2P_1(Y)Y'}{1 + e^{2Y} } + P_1(Y)Y' 
$$
defines $f$ as a meromorphic function in the plane, with no zeros. Next, set $R = f'' + b_1 f' + b_0 f$, where $b_1 = -P'/P - Y''/Y' $ and $b_0 = - (PY')^2 $.
This gives, since $h''= (Y''/Y')h'+ (Y')^2 h$,  
$$
\frac{R}f = (P-P^2) \left( (Y')^2 - \left( \frac{h'}{h} \right)^2 \right) = \frac{ (P-P^2)  (Y')^2 }{h^2} ,
$$
and so $R$ is zero-free in some
$\Omega(r_1)$.
Moreover, $S = R (P-P^2)^{-1} (Y')^{-2} $ satisfies $S/f = h^{-2}$ and $S'/S = (P-2)h'/h$.
Hence the same construction, with $P$ replaced by $P-2$,
gives rational functions  $c_j$, $d_j $ and  $e_j$  such that 
$$
\frac{S''+c_1 S'+c_0 S}{S} =
\frac{R''+d_1 R'+d_0 R}{R} = \frac{f^{(4)} + e_3 f^{(3)} + \ldots +  e_0 f}{R} = \frac{F}{R}
$$
is free of zeros in some $\Omega(r_2)$, as is $F$.

\section{Preliminaries}
\label{prelim}

\begin{lem}
\label{lemintegervalued}
Let  the function $g$ be analytic on the half-plane $H^+$ given by
${\rm Re} \, z \geq 0$, such that $g(n) \in \Z$ for all $n \in \Z \cap H^+$ and 
$ |g(z)| = o\left( 2^{|z|} \right)$ as $z \to \infty $ in $H^+$. Then $g$ is a polynomial.

Next, let $h(z) = e^{2 \pi i \alpha z} u(z)$, where $\alpha \in \R$ and $u$ is analytic on $H^+$, with 
$\log^+ |u(z)| = o( |z|)$ as $z \to \infty $ in $H^+$, and assume that $h(n) = 1$ for all large $n \in \N$. Then
$u(z) \equiv 1$ and $\alpha \in \Z$. 
\end{lem}
\textit{Proof.} The first assertion is proved in \cite{Lainteger}. To prove the second part
let $\delta_1 \in (0, \infty)$ be small: then there exist  $p \in \Z$ and $q \in \N$ such that
$$
\alpha = s + t, \quad s = \frac{p}q, \quad  |t| < \frac{\delta_1}{2 \pi q} .
$$ 
Here $t=0$ if $\alpha$ is rational, while if $\alpha \not \in \Q$ then suitable $p$ and $q$ exist by  
Dirichlet's approximation theorem \cite[p.155]{hardywright}. 
Write 
$$
z = qw, \quad F(w) = e^{2 \pi i tq w} u(qw).
$$
If $n \in \N$ is large then 
$$
1 = h(qn) = e^{2 \pi i \alpha qn} u(qn) = e^{2 \pi i (pn+ tqn)  } u(qn) = F(n) .
$$
Thus $F$ is a polynomial, by the first part, and so $F(w) \equiv 1$. Moreover, $t=0$, because otherwise there exists 
$\theta \in \{ - \pi/4, \pi/4 \}$ such that $F(re^{i \theta } ) \to 0$ as $r \to + \infty$, and so $u \equiv 1$. 
Finally, $\alpha = p/q$ must be an integer, since  $1 = h(qn+1) = \exp ( 2 \pi i p/q ) $ for large $n \in \N$.
\hfill$\Box$
\vspace{.1in}

\begin{lem}
\label{lemPL}
 Let $d_1$, $d_2$ and $\lambda $ be positive constants and let $g $  be a zero-free analytic function 
on the half-plane $Re (w) >  0$, with $\log ^+ | g(w) |  \leq d_1 + d_2 |w|^\lambda $ there.
Then for  $0 < \alpha < \pi /2 $ there exists $\mu = \mu_\alpha > 0$ such that 
$\log^+ |1/g(w)| \leq \mu_\alpha  |w|^{1+\lambda } $ as $w \to \infty $ with 
$| \arg w | \leq 
\alpha $. 
\end{lem}
\textit{Proof.} This is  standard: set
$w = (1+z)/(1-z)$ and $ g(w)=G(z)$ for $|z| < 1$. With $\rho = (1 + r)/2$ this leads to 
\begin{eqnarray*}
\log M(r, 1/G) &\leq&  \left( \frac{\rho+r}{\rho-r} \right) T( \rho, 1/G) 
\leq \left( \frac{\rho+r}{\rho-r} \right) ( \log M( \rho, G)  + O(1)) =  O( 1-r)^{-1-\lambda}
\end{eqnarray*}
as $r \to 1-$. It remains only to observe that there exists
$c_{1} = c_1(\alpha) > 0$ such that if $|w|$ is large and $ | \arg w | 
\leq \alpha < \pi /2$  then  
$(1-|z|^2)^{-1} \leq c_1 |w|$.
\hfill$\Box$
\vspace{.1in}

\begin{lem}
\label{lemb2}
 Suppose that $p$ and $q$ 
are (both formal or both locally analytic) solutions of the equations
\begin{equation}
 \label{b2}
\frac{p'}{p} = d_0 \, \frac{q'}{q} + d_1, \quad   q'' + \nu_1 q' + \nu_0  q  = 0,
\end{equation}
where the $d_j$ and $\nu_j$ are rational at infinity, and let $L$ be as in (\ref{1}). 
Then there exist  coefficients $b_j$, each  rational at infinity,
such that $p$ and $q$ satisfy 
\begin{equation}
 \label{b3}
\frac{L[p]}{p} = \sum_{j=0}^k b_j \left( \frac{q'}{q} \right)^j , \quad b_k = d_0 (d_0-1) \ldots (d_0 - k+1) .
\end{equation}

Moreover, if $d_0 \not \equiv 0, 1$ and $e_1$ and $e_0$ are rational at infinity, then 
there exist  coefficients $E_\mu$, each  rational at infinity  and depending only on the $d_j$, $e_j$ and $\nu_j$, 
such that $E_2 \not \equiv 0$ and 
\begin{equation}
 \label{b3aa}
E_2 \frac{p''}{p} + E_1 \frac{p'}{p}  + E_0  =  \left( \frac{p'}{p} \right)^2  + e_1  \frac{p'}{p}  + e_0 .
\end{equation}
\end{lem}
\textit{Proof.} Formula (\ref{b3}) follows from (\ref{b2}) and a simple induction argument, which  deliver
$$
\left( \frac{q'}{q} \right)' = -  \left( \frac{q'}{q} \right)^2 -\nu_1  \left( \frac{q'}{q} \right) - \nu_0 ,
\quad \frac{p^{(m)}}{p} = \sum_{j=0}^m b_{j,m} \left( \frac{q'}{q} \right)^j, \quad m \in \N, 
$$
with the $b_{j,m}$ 
rational at infinity and $b_{m,m} = d_0 (d_0-1) \ldots (d_0 - m+1) $.

To prove the second part, suppose that $d_0 \not \equiv 0, 1$ and 
write $P= p'/p$ and $Q= q'/q$ so that $Q = AP+B$, with $A$, $B$ rational at infinity and $A \not \equiv 0, 1$.
This yields
\begin{eqnarray*}
0 
&=& A' P + A P' + B' + A^2 P^2 + 2AB P + B^2  + \nu_1 ( AP+B ) + \nu_0 \\
&=& (A^2-A) P^2 + A (P'+P^2) + ( A'+2AB+ \nu_1A)P + B'+B^2 + \nu_1B+ \nu_0 ,
\end{eqnarray*}
and so 
\begin{eqnarray*}
(A-A^2)( P^2 + e_1 P + e_0) &=&  
A (P'+P^2) + ( A'+2AB+ \nu_1A + (A-A^2)e_1)P + \\
& & + B'+B^2 + \nu_1B+ \nu_0 + (A-A^2) e_0 .
\end{eqnarray*}
 \hfill$\Box$
\vspace{.1in}

\begin{lem}
 \label{lem3}
Suppose that $u$ and $v$ are linearly independent  (both formal or both locally analytic) solutions of an equation
$$
y''' + B_2 y'' + B_1 y' + B_0 y =0,
$$
with the $B_j$ rational at infinity, and assume that $W = W(u, v)= uv'-u'v$ is such that $W'+E_1W = 0$, where 
$E_1 $ is  rational at infinity. Then $u$ and $v$ solve an equation 
\begin{equation}
 \label{7}
y'' + E_1 y' + E_0 y = 0,
\end{equation}
where $E_0$ is also  rational at infinity. 
\end{lem}
\textit{Proof.} Since 
$u$ and $v$ are solutions of the equation $W(u, v, y) = 0$, 
it is enough to prove that $W(u',v') = E_2 W(u, v)$ with $E_2$ rational at infinity. But $W' = u v'' - u'' v = - E_1 W$
leads to 
\begin{eqnarray*}
 (-E_1'+E_1^2) W &=&  W'' = u'v'' - u''v' + uv''' - u''' v \\
&=& W(u',v') + v( B_2u'' + B_1 u' + B_0 u) - u( B_2 v'' + B_1 v' + B_0 v) \\
&=& W(u',v') - B_2 W' - B_1 W = W(u',v') + (E_1B_2 - B_1) W.
\end{eqnarray*}
\hfill$\Box$
\vspace{.1in}

\section{Asymptotics for linear differential equations} \label{wasowthm}

As in \cite{Bru} a fundamental role will be played by  formal and asymptotic expansions for solutions
of linear differential equations. 
For an equation $L[y]=0$, with $L$ as in (\ref{1}) and the $a_j$ rational
at infinity, classical results  (see  \cite[Theorem 19.1]{wasow} or \cite{Bru,jurkat})) show that 
there exist $p \in \N$ and a fundamental set of $k$ linearly independent formal solutions
\begin{equation}
 \label{formalsol}
 \widetilde h_j(z) = \exp( P_j( z^{1/p}) ) z^{\gamma_j} \sum_{\mu=0}^{n_j} U_{j,\mu}( z^{1/p} ) (\log z)^{\mu } 
\end{equation}
which satisfy the following:  
$\gamma_j $ is a complex number;  $n_j$ is a non-negative integer; the \textit{exponential part} 
$ P_j( z^{1/p})$ is a polynomial in $z^{1/p}$;
the $U_{j,\mu}(z^{1/p} ) $ are formal series in  descending integer powers of $z^{1/p}$, that is, in which at most finitely many positive powers occur;
the lead series $U_{j,n_j}$ is not the zero series. Formal solutions (not necessarily linearly independent)
with  these properties will be referred to henceforth as \textit{canonical formal solutions}.

A standard 
approach \cite{wasow} to obtaining these $  \widetilde h_j(z)  $ is to transform a solution $h$ of $L[y]=0$ into a  vector  
$\mathbf{h} = (h, h', \ldots, h^{(k-1)} ) $,
so that a fundamental solution set for $L[y]=0$ corresponds to the first row of a 
matrix solution $V(z) = U(z) z^G e^{Q(z)}$ 
of an equation $Y' = A(z) Y$, where $Q(z)$ is a diagonal matrix, its entries polynomials in $z^{1/p}$, while $G$ is a constant matrix, which may be assumed to be in Jordan
form, and $U(z)$ is a matrix with entries which are formal series in descending integer powers of $z^{1/p}$. 
Furthermore, for each $\theta \in \R$ there exists $\delta = \delta (\theta) > 0$ such that $L[y]=0$ has a fundamental set of analytic solutions 
\begin{equation}
 \label{analyticsol}
h_j(z) = \exp( P_j( z^{1/p}) ) z^{\gamma_j} \sum_{\mu=0}^{n_j} V_{j,\mu}( z^{1/p} ) (\log z)^{\mu } 
\end{equation}
on a sector $S$ given by $|z| > R_0 > 0$, $| \arg z - \theta | < \delta $, in which each $V_{j,\mu}( z^{1/p} )$ is analytic on $S$ and satisfies
$V_{j,\mu}( z^{1/p} ) \sim U_{j,\mu}( z^{1/p} )$ as $z \to \infty$ on $S$, in the sense of 
asymptotic series (see \cite[Theorem 19.1]{wasow} or  \cite{lutz}). Here
$W(z) \sim  U(z) = \sum_{m=M}^\infty U_m z^{-m/p} $ as $z \to \infty $ on  $S$
means that, for each $n \geq M$, 
$$
W(z) - \sum_{m=M}^n U_m z^{-m/p} = o( |z|^{-n/p} ) \quad \hbox{as $z \to \infty $ in $S$.}
$$

It may be assumed that the exponential parts $P_j( z^{1/p})$ have zero constant term, and this convention will
be used throughout. 
Given any exponential part $P_j( z^{1/p})$
arising for $L[y]=0$, there is always a canonical formal solution with  
exponential part $P_j( z^{1/p})$ which is free of logarithms, that is, has $n_j=0$; this holds because the matrix $G$ may be chosen to be in Jordan form. 
The  following lemma is  well known \cite{Bru,jurkat,wasow}.

\begin{lem}
 \label{wronskianlem2}
Given $k$ linearly independent canonical formal solutions of 
$L[y]=0$ with  exponential parts $q_1, \ldots, q_k$, their 
formal Wronskian has exponential part $\sum_{j=1}^k q_j$, and the exponential parts of any  
fundamental set of canonical formal solutions  of $L[y]=0$ form a permutation of the $q_j$. 
\end{lem}

\hfill$\Box$
\vspace{.1in}

For the special case of a second order equation, suppose 
that $A^*$ is rational at infinity, with $A^*(z) = (1+o(1)) c_n z^n  $ as $z \to \infty$, where $c_n \in \C^*$ and $n \geq -1$. 
Then infinity is an irregular singular point for  
\begin{equation}
 \label{A*eqn}
 w'' + A^* w = 0 ,
\end{equation}
and asymptotics are developed via Hille's 
method \cite{Hil2} as follows. The \textit{critical rays} are given by 
$\arg z = \theta^*$, where $c_n e^{i (n+2) \theta^*} $ is real and positive. 
If $0 < \beta < 2 \pi /(n+2) $ then, in a sector given by $|z| > r_1, \,  | \arg z - \theta^*| < \beta $, there exist linearly independent analytic solutions, for $j=1,2$,
\begin{equation}
 \label{hille*}
 \phi_{j}(z) = A^*(z)^{-1/4} (1+o(1)) \exp ( (-1)^j i Z ), \quad Z = \int^z A^*(t)^{1/2} \, dt = 
\frac{2 c_n^{1/2} z^{(n+2)/2}}{n+2} + \ldots .
\end{equation}
If $n=-1$ then this sector should be understood as lying on the Riemann surface of $\log z$. 
To one side of the critical ray, one of these  solutions is large and the other small,
and these roles are reversed as the critical ray is crossed. Any linear combination
$D_1 \phi_1 + D_2 \phi_2 $ with $D_1,  D_2 \in \C^*$ has a sequence of zeros tending to infinity
near the critical ray. 
Moreover, the corresponding
formal solutions, to which the $\phi_j$  are asymptotic,
may  be calculated readily from $A^*(z)$ and $Z$, with $n_j =0$ and  $p \in \{ 1, 2 \}$ in (\ref{formalsol}) (see \cite{La9} for details).

\begin{lem}
 \label{wronskianlem3}
Suppose that $F_1, \ldots, F_k$ are formal expressions, each of which is given by 
$$F_j(z)  = U_j (z) z^{\gamma_j} r_j(z) \exp( q_j(z)) ,$$
with $\gamma_j \in \C$,
$r_j$ a rational function,  $q_j$  a polynomial and $U_j(z) = 1 + O(1/z) $, in which $O(1/z)$ denotes a formal series 
in negative integer powers of $z$. Assume that none of the $r_j$ vanishes identically, 
and that $q_j - q_{j'}$ is non-constant for $j \neq j'$.
Then the formal Wronskian $W = W(F_1, \ldots, F_k) $
has an expansion
$$W(z) = \left(1 + O \left( \frac1z \right) \right) \left( \prod_{j=1}^k \left[ z^{\gamma_j } r_j(z) \exp( q_j(z) ) \right] \right)
\left( \prod_{k \geq m > n \geq 1} \left[q_m'(z) - q_{n}' (z) \right]  \right) .
$$
\end{lem}
\textit{Proof.} This is standard, and is proved by induction on $k$, using
\begin{eqnarray*}
 W(F_1, \ldots, F_k) &=& F_1^k W \left(1, \frac{F_2}{F_1} , \ldots, \frac{F_k}{F_1}  \right) =
F_1^k W \left( \left( \frac{F_2}{F_1} \right)', \ldots, \left(\frac{F_k}{F_1} \right)' \right), \\
\left( \frac{F_j}{F_1}  \right)' (z) &=& \left(1 + O \left( \frac1z \right) \right) z^{\gamma_j - \gamma_1} 
\left( \frac{r_j(z)}{r_1(z)} \right) (q_j'(z) - q_1'(z)) 
\exp( q_j(z) - q_1(z)) .
\end{eqnarray*}
\hfill$\Box$
\vspace{.1in}

\section{Beginning the proof of Theorem \ref{thm1}: Frank's method}
\label{frankmethod}

Assume that $f$ is as in the hypotheses of Theorem \ref{thm1}.
Let $f_1, \ldots, f_k$ be linearly independent 
locally analytic solutions of $L[y] = 0$.
Frank's method \cite{Fra1,FH} defines
$g, h, w_j$ and $Y$ locally by 
\begin{equation}
 \label{2}
g^k = \frac{f}{F}, \quad h = - \left( \frac{f'}{f} \right)  g, \quad w_j = f_j' g + f_j h = fg \left( \frac{f_j}f \right)' , 
\quad \frac{Y'}{Y} = - a_{k-1} . 
\end{equation}
Note that $g$ might not be meromorphic in $\Omega(r_1)$, but $g'/g$ is, and has a simple pole with residue~$1$
at every pole of $f$; moreover, at a pole of $f$ of multiplicity $m_0$, calculating the leading Laurent coefficient of $F/f$ gives 
\begin{equation}
\label{ac7}
 (g')^{-k} = (-1)^k m_0(m_0+1) \ldots (m_0+k-1)  . 
\end{equation} 

Now write locally, using Abel's identity, $W(f_1, \ldots, f_k) = cY$ and 
$$
\frac{Y }{(fg)^k}  = \frac{YF}{f^{k+1}} = \frac{ c W(f_1, \ldots, f_k, f) }{f^{k+1}}
= c W\left(\frac{f_1}f, \ldots, \frac{f_k} f, 1 \right) = c W\left( \left(\frac{f_1}f\right)', \ldots, \left(\frac{f_k}f\right)' \, \right),
$$
so that $Y = c W(w_1, \ldots, w_k)$ by (\ref{2}). 
Hence the $w_j$ are linearly independent (local) solutions of an equation $M[y] = 0$, where 
\begin{equation}
 \label{3}
M = D^k + A_{k-1}D^{k-1} +  A_{k-2}D^{k-2} + \ldots + A_0 , \quad A_{k-1} = a_{k-1} , \quad D = \frac{d \,}{dz} .
\end{equation}
Here a pivotal role is played by whether or not the differential operators $L$ and $M$ are the same, 
and Br\"uggemann's method in \cite{Bru}
depends on reducing the problem to the case $L=M$.
It will be proved in Proposition \ref{FHcase} below
that if $L=M$ then all  poles $z$ of $f$ with $|z| $ sufficiently large have the same multiplicity. 
Thus Example II in Section \ref{exa} demonstrates that 
$L$ and $M$ can indeed be different operators. 

By Frank's method, the $A_j$ are analytic in some annulus $\Omega (r_1)$ and satisfy $T(r, A_j) = S(r, f'/f) $, where 
$S(r, f'/f)$ denotes any term which is $O( \log T(r, f'/f) +  \log r ) $ 
as $r \to \infty$, possibly outside a set of 
finite measure \cite{Hay2}: see \cite[Section 2 and Lemmas A, B and 5]{FL2} for details, including
the Nevanlinna characteristic in $\Omega (r_1)$.
Denote by $\Lambda$ the field generated by the $a_j, A_j$ and their 
derivatives: then $T(r, \lambda ) 
= S(r, f'/f)$ for all $\lambda \in \Lambda$. 

To simplify the subsequent calculations  it is convenient to write 
\begin{equation}
 \label{gf1}
- k \, \frac{X'}{X}  =   a_{k-1} =   A_{k-1}, \quad p = \frac{f}X , \quad  p_j = \frac{f_j}X , \quad q  = - \left( \frac{p'}{p} \right) g, \quad 
t_j = p_j' g + p_j q = \frac{ w_j }X .
\end{equation}
It is then well known that there exist equations 
\begin{equation}
 \label{gf2}
y^{(k)} + c_{k-2} y^{(k-2)} + \ldots + c_0 y = 0,
\end{equation}
and 
\begin{equation}
 \label{gf3}
y^{(k)} + C_{k-2} y^{(k-2)} + \ldots + C_0 y = 0,
\end{equation}
in which the $c_j$ and $C_j$ all belong to $\Lambda$, and 
$c_{k-2} = C_{k-2}$ if and only if $a_{k-2} = A_{k-2}$, 
such that $L[Xy]=0$ if and only if $ y $ solves (\ref{gf2}), and 
$M[Xy]=0$ if and only if $ y $ solves (\ref{gf3}). In particular,
the $p_j$ and  $t_j$ are linearly independent local  solutions of (\ref{gf2}) and  (\ref{gf3}) respectively. 
The following lemma \cite{FH} is key to Frank's method:
see  also \cite[Lemma C]{FL2}. 

\begin{lem}[\cite{FH}]
\label{franklem1}
Let $G$, $\Phi$, $p_1, \ldots , p_k$,
$ c_0, \ldots , c_{k-2}$ and $C_0, \ldots , C_{k-2}$  be analytic functions on a plane domain $U$, such that 
$p_1, \ldots, p_k$ are linearly independent solutions of (\ref{gf2}). 
Then the functions 
$ p_1' G + p_1 \Phi, \ldots ,  p_k' G + p_k \Phi$ are 
solutions in $U$ of the equation $(\ref{gf3})$
if and only if, with the notation $C_k = 1$ and
$c_{k-1} = C_{k-1}= c_{-1} = C_{-1} = 0$
and
$$
M_{k,\mu}[w] =
\sum_{m=\mu}^k \frac{m!}{\mu!(m-\mu)!}  C_m w^{(m-\mu)} \quad (0 \leq \mu \leq k),
\quad M_{k,-1}[w] = 0,
$$
the functions $G$ and $\Phi$ satisfy, for $0 \leq \mu \leq k-1$,
\begin{equation}
M_{k,\mu} [\Phi] - c_\mu \Phi = -
M_{k,\mu-1}[G]  + c_\mu M_{k,k-1}[G] +
(  c_\mu' + c_{\mu-1} ) G .
\label{9818}
\end{equation}
\end{lem}
\hfill$\Box$
\vspace{.1in}

Because the proof of Lemma \ref{franklem1} is based on purely formal calculations,
an analogous statement holds linking formal solutions of (\ref{gf2}), (\ref{gf3}) and (\ref{9818}). 
Since the coefficient of $\Phi$ in $M_{k,\mu} [\Phi] - c_\mu \Phi$ is $c_{0,\mu} = D_\mu = C_\mu - c_\mu$,
the equations (\ref{9818}) may be written in the form 
\begin{equation}
T_\mu [G] = S_\mu [\Phi] =
\sum_{j=0}^{k - \mu } c_{j, \mu} \Phi^{(j)} ,
\quad c_{0,\mu} = D_\mu = C_\mu - c_\mu  ,\quad 0 \leq \mu \leq k-1, 
\label{9821}
\end{equation}
in which $T_\mu$ and $S_\mu$ are linear differential
operators with coefficients in $\Lambda $.
In particular these equations are satisfied by $G=g$, $\Phi=q$. 
Taking  $\mu = k-1$ in (\ref{9818}) produces 
 \begin{equation}
-\Phi' = - U[G] = 
\frac{(k-1)G''}2 + \frac{D_{k-2} G}k .
\label{9821X}
\end{equation}

Now $\mu = k-2$ and (\ref{9821X}) give (as in \cite[p.162]{GL} or \cite[Lemma 8, pp.307-8]{Latsuji})
\begin{equation}
 \label{9821Y} 
D_{k-2}  \Phi = \frac{k(k^2-1)}{12} G''' + G' \left( -  \frac{(k+1)D_{k-2}}2   + 2 C_{k-2} \right) +  G \left( \frac{k-1}2 D_{k-2}' + c_{k-2}' - D_{k-3} \right) .
\end{equation}
Next, combining (\ref{9821X}) with (\ref{9818}) for $\mu = k-3$ yields, with $d_j$ denoting elements of $\Lambda$, 
\begin{equation}
 \label{9821Z} 
\frac{2 D_{k-3}}{k-2} \Phi = \frac{k(k^2-1)}{12} G^{(4)} + G'' \left( \frac{(k-1)D_{k-2}}3  + 2 C_{k-2} \right) + d_1 G' + d_2 G. 
\end{equation}
Note that (\ref{9821Z}) holds even if $k=3$, in which case $M_{k,k-4} [G] = 0$. 
Differentiating (\ref{9821Y}) and using (\ref{9821X}) and (\ref{9821Z}) leads to
\begin{equation}
D^* \Phi = \left( \frac{2 D_{k-3}}{k-2} - D_{k-2}' \right) \Phi =
\frac{(k+2) D_{k-2}}3 \,  G'' + d_3 G' + d_4 G. 
\label{D*eqnQ}
\end{equation}

\begin{lem}
 \label{lemgDE}
There exists a non-trivial homogeneous linear differential equation $N_1[y]=0$,  of order at most~$3$
and with coefficients in $\Lambda$,  
with the property that if the pair $\{ G, \Phi \}$ solves the system (\ref{9821}) then
$G$ solves $N_1[y]=0$. 
\end{lem}
\textit{Proof.} If $D_{k-2} \equiv 0$ this is clear from (\ref{9821Y}), so assume that 
$D_{k-2} \not \equiv 0$. 
If $D^*$ vanishes identically in (\ref{D*eqnQ}) then  a second order equation arises for $G$, while otherwise combining (\ref{D*eqnQ})  with 
(\ref{9821Y}) yields a 
third order equation. 
\hfill$\Box$
\vspace{.1in}

Consider now two cases. 
\\\\
{\bf Case 1.} Assume that $c_{0,\mu} = C_\mu - c_\mu \equiv 0$ for
$0 \leq \mu \leq k-1$ in (\ref{9821}). 
This is equivalent to the 
equations (\ref{gf2}) and (\ref{gf3}) being the same,
and hence equivalent to the 
operators $L$ and $M$ being identical. 
In this case, 
$t_j = p_j ' g + p_j q $ is a solution of (\ref{gf2})
for $1 \leq j \leq k$.
Since $p_1$ and $p_2$ are linearly independent,
$p_1 p_2' - p_1'p_2 $ does not vanish identically and so (\ref{gf1}) yields
$$
H_1 = \frac{f'}f + \frac{a_{k-1}}k = \frac{p'}p = - \frac{q}g = \frac{ p_1' t_2 - p_2 ' t_1 }{p_1 t_2 - p_2 t_1 }.
$$
For $\theta \in  \R$ and $\kappa \in \C$ the number of distinct zeros of $\kappa - H_1$ in
$r_1 +1 \leq  |z| \leq r, | \arg z - \theta | \leq \pi/4 $, is at most the number of zeros of 
$p_1' t_2 - p_2 ' t_1 - \kappa ( p_1 t_2 - p_2 t_1 )$ there, which is
bounded by a power of $r$ as $r \to \infty$, by  \cite[Lemma 2]{FL2} or standard sectorial methods.
Hence  $f'/f$ has finite order of growth, by the second fundamental theorem, and  every $\lambda \in \Lambda$ is rational at infinity. 
\\\\
{\bf Case 2.} Assume that the coefficient of
$\Phi$ in at least one of the $S_\mu $ in (\ref{9821}) is not 
identically zero, this being equivalent to $L \neq M$. 

\bigskip

Let $\nu $ be the largest integer with $0 \leq \nu \leq k-1$
such that $c_{0, \nu} \not \equiv 0$. Then every pair $\{ G, \Phi \}$ satisfying the system (\ref{9821}) 
(including $\{ g, q \}$) 
has 
\begin{equation}
\Phi = (c_{0, \nu})^{-1} \left( 
T_\nu[G] - 
\sum_{j=1}^{k-\nu} c_{j, \nu} \frac{d^{j-1}}{dz^{j-1}}(U[G])\right) = T^*[G] ,
\label{helim}
\end{equation}
by (\ref{9821X}). Observe that the operator $T^*$ has order at least $1$,
and in particular is not the zero operator, since otherwise (\ref{gf1}) leads to 
$ (-p'/p) g = q = \eta_1 g$, with $\eta_1 \in \Lambda $, so that $f'/f \in \Lambda$, and hence
$f'/f$ is rational at infinity, contrary to assumption.
It follows from (\ref{9821}), (\ref{9821X}) and (\ref{helim})
that if $\{ G, \Phi \}$ solves  (\ref{9821}) then $G$  solves
the system
\begin{equation}
U[G] = \frac{d}{dz}(T^*[G]), \quad 
S_\mu ( T^*[G] ) = T_\mu [G] , \quad 0 \leq \mu \leq k-2 ,
\label{9825}
\end{equation}
as does, in particular, $g$.
Conversely, if $G$ solves the system (\ref{9825})
(in the analytic or formal sense), then (\ref{9821}) is satisfied by setting $\Phi = T^*[G]$.
This system (\ref{9825}) cannot be trivial, because otherwise (\ref{9821}) holds 
with  $\Phi = T^*[G]$ and an arbitrary choice of $G$, which would then have to solve the equation $N_1[y]=0$ of Lemma \ref{lemgDE}.
A standard reduction procedure \cite[p.126]{Ince} now generates a non-trivial homogeneous linear differential
equation $N[y]=0$, with coefficients in the  field $\Lambda$, 
whose (analytic or formal) solution space coincides with that of the system 
(\ref{9825}). Here every solution $G$
of $N[y]=0$ is such that the pair $\{ G, T^*[G] \}$ solves (\ref{9821}),
and so $G$  solves  $N_1[y]=0$, from which it follows that
$N$ has order at most $3$.

Suppose that $N$ has order $1$: then $g'/g \in \Lambda$, and so $p'/p$ and $f'/f$ belong to $\Lambda$, by (\ref{gf1}) and (\ref{helim}),
so that $f'/f$ is rational at infinity, contrary to assumption. 
Thus $N$ has order at least $2$, but at most $3$, and the system (\ref{9825}) has a solution $G$ with $G/g$ non-constant.
By an argument from \cite{FH}
(see \cite[Proof of Theorem 3, Case 1B]{FL2} for details), 
$p'/p$ has a representation as a rational function in the $p_j$ and their derivatives. 
The same sectorial argument as used in Case 1 shows that  $f'/f$  has finite order of growth, as has $g^k = f/F$, and
all members of the field $\Lambda$ are rational at infinity. 

Hence the fact that $N$ has order at most $3$ gives  an operator $V_2$, having order at most $2$ and 
coefficients which are
rational at infinity, with the following property. Every solution 
$G$ of (\ref{9825}) has $T^*[G] = V_2[G]$, so that the pair $\{ G, V_2[G] \}$ solves (\ref{9821}),
and $p_j' G + p_j V_2[G] $ solves (\ref{gf3}) for $j=1, \ldots , k$, by Lemma~\ref{franklem1}.
Moreover, (\ref{helim}) gives $q = T^*[g] = V_2[g]$. With $V = V_2 + a_{k-1}/k$, this implies using (\ref{gf1}) that each 
$f_j' G + f_j V[G]$ solves $M[y]=0$, and 
$-f'/f = h/g = V[g]/g$.
It now follows,
using the Wiman-Valiron theory \cite{Hay5} and estimates for logarithmic derivatives \cite{Gun}
applied to $g$, 
that $f$ has finite order. It also follows that
$f$ has   an unbounded sequence of poles, since otherwise $f'/f$ is rational at infinity. 
The following key lemma has thus been proved. 

\begin{lem}
 \label{lemgeqn}
 With the hypotheses of Theorem \ref{thm1}, the function $f'/f$ has finite order, 
 all elements of the field $\Lambda$ are rational at infinity, and
$h = - (f'/f)g $ satisfies 
\begin{equation}
 \label{5a}
- h' = \left( \frac{k-1}2 \right) \, g'' - \frac{a_{k-1}}k g' + \frac{D_{k-2}-a_{k-1}'}k g, \quad D_j = C_j - c_j .
\end{equation}

Furthermore, if the operators $L$ and $M$ are not the same then the following additional conclusions hold.
The function $g$ solves a homogeneous linear differential equation
$N[y] = 0$, of order $2$ or $3$, with coefficients which are rational at infinity. 
Moreover, $f$ has finite order and an unbounded sequence of poles,
and there exist functions $\alpha, \beta, \gamma$, 
all rational at infinity,
such that
\begin{equation}
 \label{5}
h = - \left( \frac{f'}f \right)  g = V[g], \quad V = \alpha D^2 + \beta D + \gamma .
\end{equation}
Finally, if $G$ is a locally analytic solution of $N[y]=0$, and $K$  is a locally analytic solution of $L[y]=0$, then 
$K' G + K V[G]$ is a (possibly trivial) solution of $M[y]=0$. 
\end{lem}

\hfill$\Box$
\vspace{.1in}

Here (\ref{5a}) follows from (\ref{gf1}) and (\ref{9821X}). 
The last assertion of Lemma \ref{lemgeqn}
also holds for  formal solutions $G$ and $K$ of $N[y]=0$ and $L[y]=0$ respectively.

\section{The first special case }

\begin{prop}
 \label{FHcase}
With the hypotheses of Theorem \ref{thm1}, suppose in addition that $c_{k-2} = C_{k-2}$
in (\ref{gf2}) and (\ref{gf3}), which holds 
if and only if $a_{k-2} = A_{k-2}$ in $L$ and $M$, and certainly holds
if the operators $L$ and $M$ are the same. Then $f$ satisfies conclusion (i) of Theorem \ref{thm1}. 
\end{prop}
\textit{Proof.} 
The approach here is essentially due to Frank and Hellerstein \cite{FH}. Since $D_{k-2} =0$ in (\ref{5a}),
integration gives a constant $d$ such that  
\begin{equation}
\frac{f'}f = - \frac{h}{g} = \frac{(k-1)g'}{2g } + \frac{d}g - \frac{a_{k-1} }k   . 
\label{repr}
\end{equation}
If $d = 0$ then comparing residues shows that $f$
has no poles in some $\Omega(r_2)$ and $F/f$, which has finite order by Lemma \ref{lemgeqn}, satisfies  
$g^{-k} = F/f = R_1 e^{P_1}$ with $R_1$  rational at infinity and $P_1$
a polynomial, so that $f'/f$ is rational at infinity, by (\ref{repr}), contrary to assumption. 

Assume henceforth that $d \not = 0$ in (\ref{repr}), which makes  $g$ 
meromorphic of finite order
in $\Omega(r_1)$.
Suppose that $f$ has no poles in some $\Omega(r_2)$. Then $g$ has no
zeros and poles there and $g = R_2 e^{P_2}$ 
in (\ref{repr}), with $R_2$  rational at infinity and $P_2$
a polynomial. This gives, since $f'/f$ is not rational at infinity,  
\begin{equation}
 \label{Brurep6}
\frac{f'}{f} = R + Se^P , \quad SP' \not \equiv 0 ,
\end{equation}
where  $R$ and $ S$ are rational at infinity and $P$ is a polynomial. 
It follows using \cite[Lemma 3.5]{Hay2} that 
\begin{equation}
 \label{Brurep3}
\frac1{g^k} =
\frac{F}{f} = S^k e^{kP} + e^{(k-1)P} \left( k S^{k-1} R + a_{k-1}S^{k-1} + \frac{k(k-1)}2 S^{k-2} (S' + P'S) \right)  + \ldots  . 
\end{equation}
Since $F/f $ has neither zeros nor poles in $\Omega(r_2)$, the coefficient of $e^{(k-1)P}$ must vanish identically,
leading to the first equation of  (\ref{Brurep}), with $H'=Se^P$, and to $F/f = S^k e^{kP} =(H')^{k}  $. 
Here $H''/H'$ does not vanish at infinity, because $P'$ does not. 

Suppose next that $f$ has an unbounded sequence of poles.
At a pole $z$
of $f$, with $|z|$ large and with multiplicity
$m$, equations (\ref{ac7}) and (\ref{repr})   deliver 
$$
\frac1{d^{k}} = \chi (m) = \frac{ m(m+1) \ldots (m + k-1 )}{(m + (k-1)/2 )^{k}}, 
$$
so that $d^k $ must be real and greater than 1, by the arithmetic-geometric mean
inequality.
A further application of the same
inequality to $\chi ' / \chi $ shows that all poles $z$ of $f$ with $|z|$ sufficiently large 
have fixed multiplicity $m$.
Set $T_1 = f'/f$.
Since $g^k$ and $T_1$ have
finite order, standard estimates \cite{Gun} give
$M_1 > 0$ such that $T_1^{(j)}(z)/T_1(z) = O( |z|^{M_1} )$ and 
$g^{(j)}(z)/g(z) = O( |z|^{M_1} )$, for $|z|$ outside a set $F_0$ of finite
measure and
$1 \leq j \leq k$. If $|z| \not \in F_0$ and 
$\log^+ |T_1(z)| / \log |z| $ is sufficiently large this leads, using (\ref{repr}) and \cite[Lemma 3.5]{Hay2},  to 
$$\frac1{g(z)^{k}} = \frac{F(z)}{f(z)} = T_1(z)^k + \ldots = (1+o(1)) T_1(z)^k = (1+o(1))  \frac{d^k}{ g(z)^{k}},$$
which is  a contradiction since $d^k > 1$. Thus 
$\log^+ | f'(z)/f(z) | = O( \log |z| )$ for $|z|$ outside a set of finite
measure and applying the Wiman-Valiron theory 
\cite{Hay5} 
to $1/f$ shows that $f$ has finite order. Furthermore, 
since $f$ and $g$ have finite order and all poles $z$ of $f$ with $|z|$ sufficiently large 
have fixed multiplicity $m$, the
function $G_0 = f'/f + m g'/g $ is  rational at infinity. 
Substituting $f'/f = - m g'/g + G_0  $ into 
(\ref{repr})
produces a first order linear differential equation
for $g$ of form 
$$
g' + \delta g = d_0 ,
$$
with $d_0 \in \C$ and  $\delta$ rational at infinity, and with $\delta(\infty) \neq 0$, because $f/F = g^k$ has an essential singularity at infinity.  
This equation may be solved to give
$g = d_0 H/H'$, where $H''/H'=\delta$ and $H(z) \neq \infty$ and $H'(z) \neq 0$ for large $z$ 
in a sector 
containing an unbounded sequence of poles of $f$. 
It follows, using (\ref{repr}) again, that 
\begin{equation}
 \label{Brurep4}
\frac{g'}g = \frac{H'}H - \frac{H''}{H'} ,  \quad
\frac{f'}{f} = - \frac{a_{k-1}}k - \left( \frac{k-1}2 \right)  \frac{H''}{H'} +  d_1 \frac{H'}{H} , \quad d_1 \in \C.
\end{equation}
Now comparing residues shows that $d_1 = - m$ in (\ref{Brurep4}), giving the second equation of (\ref{Brurep}). 

To determine the solutions of $L[y]=0$, write 
$$
\phi = (H')^{-k} e^H , \quad \psi = (H')^{-k} H^{-m} , \quad
\delta =\frac{H''}{H'}  , 
\quad M_k = (D + \delta ) \ldots (D + k \delta ) , 
\quad D = \frac{d\,}{dz} .
$$
Then it is easy to verify that 
\begin{equation}
 \label{Brurep5}
\Phi = M_k[\phi] =  e^H, \quad  \Psi = M_k[\psi] =  c H^{-m-k}  , \quad M_k \left[ (H')^{-k} P_{k-1}(H) \right]=0, 
\end{equation}
where $P_{k-1}$ denotes any polynomial of degree at most $k-1$. 
In fact, the action of the differential operator $M_k$ on $\phi$, $\psi$ and $(H')^{-k} P_{k-1}(H)$ amounts to $k$ times differentiating with respect to $H$
the terms $e^H$,  $H^{-m}$ and $P_{k-1}(H)$. Define $Z$ locally by
\begin{equation*}
 \label{Zdef1}
\frac{Z'}Z = - \frac{a_{k-1}}k + \left( \frac{k+1}2 \right)  \frac{H''}{H'} .
\end{equation*}
Then a standard change of variables gives $L_k = D^k + \ldots + \widetilde A_1 D  + \widetilde A_0 $, 
with coefficients which are readily computable and rational at infinity,
such that 
$L_k[Zy] = Z M_k[y]$, and the last equation of
(\ref{Brurep5}) shows that $L_k[w]=0$ has linearly independent solutions $y_j$ given locally by (\ref{Brurep2}).

The next step is to show that $L_k = L$. When $f$ has no poles in some $\Omega(r_2)$, combining the first equation of (\ref{Brurep}) with 
(\ref{Brurep3}) and the remarks immediately following it yields 
$$
Z \phi = cf, \quad
\frac{L[f]}f = \frac{F}f = S^k e^{kP} = (H')^k = \frac{\Phi}{\phi}= \frac{M_k[\phi]}\phi = \frac{Z M_k[\phi]}{cf}  = \frac{L_k[cf]}{cf} = \frac{L_k[f]}f.
$$
Thus the operators $L$ and $L_k$ must agree: otherwise $f$ satisfies a homogeneous linear differential
equation with coefficients which are rational at infinity, and so has finite order, 
contradicting (\ref{Brurep6}). 
On the other hand, when $f$ has an unbounded sequence of poles, (\ref{Brurep}) and (\ref{Brurep4}) lead to
$$
Z \psi = cf, \quad
\frac{L[f]}f = \frac{F}f = \frac1{g^k} =
c \left( \frac{H'}{H}\right)^k  =    \frac{c \Psi}{\psi} =  \frac{c M_k[\psi]}{\psi} =  \frac{c Z M_k[\psi]}{f} = \frac{c L_k[f]}f.
$$
Again the operators $L$ and $c L_k$ must agree, and $c$ must be $1$,  because otherwise $f$ cannot have an unbounded  sequence of poles.
Thus, in both cases, the $y_j$ solve $L[y]=0$. 
Next, using (\ref{Brurep}) and (\ref{Brurep2}) shows, after  multiplying $y_2$ by a constant if necessary,  that 
$$
\frac{f'}{f} = \frac{y_1'}{y_1} +  \, \left( \frac{y_2}{y_1} \right)' \quad \hbox{or} \quad 
\frac{f'}{f} = \frac{y_1'}{y_1} -m \left( \frac{y_2'}{y_2} - \frac{y_1'}{y_1} \right) .
$$
This gives $f = c y_1 \exp(y_2/y_1)$ or $f = c y_1^{m+1} y_2^{-m}$ as asserted. 

Finally, set $\widetilde M_2 = (D + (k-1) \delta ) ( D + k \delta )$. There exists an operator 
$\widetilde L_2 = D^2 + \widetilde a_1 D + \widetilde a_0$, with coefficients which are rational at infinity, such that 
$\widetilde L_2[Zy] = Z \widetilde M_2 [y]$ and 
$$
\widetilde L_2[Z \phi ] = Z  \widetilde M_2 [\phi] = Z (H')^{2-k} e^H, \quad \widetilde L_2[Z \psi ] = Z  \widetilde M_2 [\psi] = c Z(H')^{2-k} H^{-m-2} .
$$
Since $f$ equals $c Z \phi$ or $c Z \psi$, there exists $r_3 > 0$ such that $\widetilde L_2 [f] $ has no zeros in $\Omega(r_3)$, and Proposition 
\ref{FHcase} is proved. 

 \hfill$\Box$
\vspace{.1in}

\section{Annihilators } 

The remainder of the proof of Theorem \ref{thm1} focuses on the case where the operators $L$ and $M$ differ.
In this case Lemma \ref{lemgeqn} ensures that if $\phi$ is a non-trivial solution of $L[y] = 0$, and 
$\psi$ is a non-trivial solution of $N[y] = 0$, then 
$\chi= \phi' \, \psi + \phi \, V[ \psi ] $ solves $M[y]=0$.
Here $\chi$ may
vanish identically, in which case  
$\psi$ will be said to annihilate $\phi$, and vice versa. This notion makes sense 
when  $\phi$ and 
$\psi$ are both analytic solutions, and also when they are both formal solutions. 
The terminology in this section is as in Section \ref{wasowthm}, and the convention that 
exponential parts  have zero constant term still applies. The following variant of an auxiliary result from \cite{Bru} is key to 
the proof of Theorem \ref{thm1}.

\begin{lem}[\cite{Bru}]
 \label{lemexppart}
 Assume that $L \neq M$ and take a canonical formal solution $G $ of  $N[y] =0$ which is free of logarithms and has
exponential part $\kappa$. In addition, take a fundamental set of canonical formal solutions $f_1, \ldots, f_k$ of $L[y]=0$, such that $f_j$ has 
exponential part $q_j$,
and a fundamental set of canonical formal solutions $w_1, \ldots, w_k$ of $M[y]=0$, where
$w_j$ has exponential part $s_j$. 
Then the following conclusions hold. \\
(i) Each $W_j = f_j' G + f_j V[G] $ is either 
identically zero or a canonical formal solution of $M[y]=0$  with exponential part 
$ q_j + \kappa $.\\
(ii) There exists $\lambda = \lambda (G) \in \{ 1, \ldots, k \}$ such that the collection $s_1, \ldots, s_k$ consists of 
\begin{equation}
 \label{sjlist1}
 q_j + \kappa \quad (j \neq \lambda), \quad q_\lambda - (k-1) \kappa .
\end{equation}
(iii) If the $W_j$ are linearly dependent,  
then $G$ annihilates a canonical formal solution $g_1$ of $L[y]=0$ with exponential part $q_\lambda$, and every formal solution of $L[y]=0$ which is annihilated
by $G$ is a constant multiple of $g_1$.\\
(iv) If  $\kappa $ is not identically  zero, then the $W_j$ are linearly dependent. 
\end{lem}
\textit{Proof.} 
Conclusion (i) follows immediately from Lemma \ref{lemgeqn}. Next, 
Lemma \ref{wronskianlem2} and Abel's identity give, 
since $a_{k-1}(\infty) = A_{k-1}(\infty)=0$, 
\begin{equation}
 \label{abel}
\sum_{j=1}^k q_j = \sum_{j=1}^k s_j = 0. 
\end{equation}

Suppose first that the $W_j$ are linearly independent. Then (i), 
(\ref{abel}) and Lemma \ref{wronskianlem2}  yield
$$
0 = \sum_{j=1}^k s_j =  \sum_{j=1}^k (q_j + \kappa ) = k \kappa  ,
$$
which implies that $\kappa = 0$ and that 
$\{ s_1, \ldots, s_k \} = \{ q_1, \ldots, q_k \}$, again by Lemma \ref{wronskianlem2}. This proves conclusion (iv), and that (\ref{sjlist1}) applies when the
$W_j$ are linearly independent.

Now suppose that the $W_j$ are linearly dependent: then $G$ annihilates a non-trivial solution $g_1 $ of $L[y]=0$. 
It may be assumed that the exponential parts and formal series appearing in $G$ and the $f_j$ and $w_j$   all involve integer powers of 
$z^{1/p}$, for some fixed $p \in \N$. 
Because $G$ is free of logarithms, (\ref{5}) implies that $V[G]/G$ is a formal series in descending powers of $z^{1/p}$, and therefore so is 
$g_1'/g_1$. Thus $g_1$ is a canonical formal solution of $L[y]=0$, and by Lemma~\ref{wronskianlem2}
it may be assumed   that $g_1=f_1$; moreover, every formal solution $g_2$
of $L[y]=0$ which is annihilated
by $G$ has $W(g_1, g_2) = 0$, so that $g_2$ is a constant multiple of $g_1$. This proves (iii).

Now set $U_j = f_j' G + f_j V[G]$. Then $U_1 \equiv 0$, but $U_2, \ldots , U_k$ are linearly independent, and 
$M[y]=0$ has a fundamental set $\{ U^*, U_2, \ldots, U_k \}$ of canonical formal solutions, 
with exponential parts $s^*, q_2 + \kappa, \ldots, q_k + \kappa$ respectively. Using (\ref{abel}) twice, as well as Lemma \ref{wronskianlem2}, shows that
these exponential parts have sum $0$ and
$s^* = q_1 - (k-1)\kappa$, which leads to (\ref{sjlist1}).
\hfill$\Box$
\vspace{.1in}

The following lemma, in which transcendentally fast means faster than any power of $z$, gives a sufficient condition for an analytic solution of $N[y]=0$ to annihilate
a solution of $L[y]=0$.

\begin{lem}
 \label{lem1}
Suppose that $L \neq M$. 
Then $g(z)$ cannot tend to $0$ transcendentally fast as $z \to \infty$ in a sector, and the equation $N[y]=0$ cannot have a fundamental set of 
canonical formal solutions with the same exponential part. 
Moreover, if $G$ is a non-trivial analytic
solution of $N[y] = 0$ and $G(z)$ 
tends to $0$ transcendentally fast as $z \to \infty$  in a sector $S$, then $G$ annihilates a non-trivial analytic solution of $L[y] = 0$. 
\end{lem}
\textit{Proof.} If $g$ tends to zero transcendentally fast on a sector, then $F/f = g^{-k}$ tends to infinity 
transcendentally fast there; since $f$ has finite order by Lemma \ref{lemgeqn}, this contradicts
standard estimates \cite{Gun} for logarithmic derivatives $f^{(j)}/f$. 

Next, if $N[y]=0$ has a fundamental set of 
canonical formal solutions with the same exponential part $\kappa$, then $\kappa$ is a polynomial in $z$,
by Lemma \ref{wronskianlem2} and Abel's identity. 
Here $\kappa$ cannot be the zero polynomial, because $g^k$ is transcendental, and so there exists a sector on which every solution of $N[y]=0$,
including $g$, tends to zero transcendentally fast, which is a contradiction.

Assume now that $G$ is a non-trivial analytic
solution of $N[y] = 0$ which 
tends to $0$ transcendentally fast in a sector $S$, but annihilates no 
non-trivial solution of $L[y] = 0$. Then there exist $k$ solutions $f_j$ of $L[y]=0$ such that the
$f_j'G + f_jV[G]$ are linearly independent solutions of $M[y] = 0$ on $S$. 
Because $N[y]=0$ has order at most $3$ and at least two distinct exponential parts, 
the asymptotics in Section \ref{wasowthm} give rise to a subsector $S^*$ of $S$
on which $G(z) \neq 0$ and $G^{(j)}(z)/G(z) = O( |z|^q )$ as $z \to \infty$, for some $q \in \N$ and all $j \in \{1, \ldots, k \}$. 
This is clear if one solution $h_j$ as in (\ref{analyticsol}) dominates the others on a subsector,
and so evidently holds unless there are two solutions $h_j$ as in (\ref{analyticsol}), with the same exponential part, for which the
powers $\gamma_j$
differ by $\delta \in i \R \setminus \{ 0 \}$; but in this case, for any given $A \in \C^*$,  a subsector 
may be chosen on which $\log |z^\delta - A| $ is bounded. 
Define 
functions $Y$, $\phi$ and $\Phi$ on $S^*$ by 
\begin{equation}
 \label{6}
 \frac{Y'}Y = - a_{k-1} = - A_{k-1}, \quad 
- \frac{\phi'}{\phi} = \frac{V[G]}{G} = \alpha \frac{G''}{G} + \beta \frac{G'}{G} + \gamma , \quad \Phi = L[ \phi ]  .
\end{equation}
It follows that
\begin{eqnarray*}
 c Y &=& W( f_1'G + f_1 V[G], \ldots , f_k'G + f_k V[G] ) \\
&=& W( f_1'G - f_1 (\phi'/\phi) G , \ldots , f_k'G - f_k (\phi'/\phi) G  ) = (\phi G)^k W( (f_1/\phi)', \ldots ,(f_k/\phi)') \\ &=& 
(\phi G)^k W(1, f_1/\phi, \ldots ,f_k/\phi) 
= \phi^{-1} G^k W(\phi, f_1, \ldots , f_k) .
\end{eqnarray*}
This delivers in turn
$$
\frac{\Phi}{\phi} = \frac{L[\phi]}{\phi} =  \frac{cW(f_1, \ldots, f_k, \phi)}{Y \phi} = \frac{c}{ G^k} ,
$$
so that $\Phi(z)/\phi(z)$ tends to infinity transcendentally fast
in the sector $S^*$. But (\ref{6}) implies that there exist  $q' , q''\in \N$ with $\phi'(z)/\phi(z)= O( |z|^{q'} )$ as $z \to \infty$ in $S^{*}$, and hence
$\Phi(z)/\phi(z)= O( |z|^{q''} )$ as $z \to \infty$
on a  subsector of $S^*$,   a contradiction.
\hfill$\Box$
\vspace{.1in}

\section{The second special case}\label{second}

\begin{prop}
 \label{frankcase}
With the hypotheses of Theorem \ref{thm1}, suppose in addition that  $L \neq M$ and that 
there exist $E \in \N$ and a function $R$ which is rational
at infinity such that all poles 
$z$ of  $\mathfrak{f}(z) = f\left(z^E\right)$ with $|z| $ sufficiently large
have multiplicity $R(z)$. 
Then $f$ satisfies at least one of conclusions (i) and (ii) of Theorem \ref{thm1}.  
\end{prop}
The proof of Proposition \ref{frankcase} will occupy the remainder of this section. 
Observe first that $f$ has finite order and an unbounded sequence of poles, by Lemma \ref{lemgeqn}.
Next, 
it may be assumed that $E=1$. To see this, let $\omega = \exp( 2 \pi i /E)$ 
and let $z_0$ be large and a pole of $f$ of multiplicity $m_0$. Let 
$w_0^E = z_0$. Then $w_0 $ is a pole of $\mathfrak{f} $ of multiplicity $m_0 = R(w_0)$. 
This is true for all $E$ choices of $w_0$ and so $R(z) = R( \omega z)$ for all large $z$, which
gives $R(z) = S\left(z^E\right)$ for some function $S$ which is rational at infinity. Thus  the multiplicity
$m_0$ of the pole of $f$ at $z_0$ satisfies $m_0 = R(w_0) = S(w_0^E) = S(z_0)$. Assume for the remainder
of this section that $E=1$. 

\begin{lem}
 \label{d0lem}
There exist functions $d_0$, $d_1$, both 
rational at infinity, such that 
$f$ and $g$ satisfy 
\begin{equation}
 \label{b0*}
\frac{f'}{f} = d_0 \, \frac{g'}{g} + d_1 .
\end{equation}
Moreover, $d_0$ either has $d_0(\infty) = \infty$ or is constant and equal to a negative integer. 
\end{lem}
\textit{Proof.} Let $d_0 = -R$. By the remark following (\ref{2}),
there exists $r_0 > 0$ such that 
$f'/f - d_0 g'/g $ has no poles in $\Omega(r_0)$, and so is rational at infinity since $g^k$ and $f$ have finite order. 
The last assertion follows  from the fact that $f$ has an unbounded sequence of poles. 
\hfill$\Box$
\vspace{.1in}

\begin{lem}
 \label{d0lem2}
There exist functions $\nu_1$, $\nu_0$,  both rational at infinity,  such that 
$g$ satisfies 
(\ref{b1}). 
\end{lem}
\textit{Proof.}  The equation (\ref{b0*}) yields, using (\ref{2}),  
$$
-h = \left( \frac{f'}f \right) g =  d_0 g' + d_1 g, \quad - h' =  d_0 g'' + (d_0' + d_1) g' + d_1' g .
$$
Combining this with  (\ref{5a}) gives 
\begin{equation}
 \label{newgeqn}
 0 = \left(d_0 - \frac{k-1}2\right) g'' + \left(d_0'+d_1 + \frac{a_{k-1}}k \right) g' + 
 \left(d_1' + \frac{a_{k-1}' + c_{k-2} - C_{k-2}}k \right) g 
\end{equation}
and an equation (\ref{b1}), as asserted,  since $d_0 - (k-1)/2 \not \equiv 0$ by Lemma~\ref{d0lem}.

\hfill$\Box$
\vspace{.1in}

From (\ref{b1}), (\ref{5})  and (\ref{b0*}) it follows that
$$
-d_0 g' -d_1 g  = - \left( \frac{f'}f \right) g = h =
V[g] = (\beta - \alpha \nu_1) g' +  (\gamma - \alpha  \nu_0 ) g
$$
and so, since $g^k$ has an unbounded sequence of zeros,
\begin{equation}
 \label{d2def}
 -d_0  = \beta - \alpha \nu_1, \quad -d_1 = \gamma - \alpha \nu_0 .
\end{equation}
In the next lemma the convention that exponential parts have zero constant term is retained. 

\begin{lem}
\label{lemb3*}
There exists an equation (\ref{A*eqn}),
with $A^*$ rational at infinity, such that $y U^{-1/2} $ solves (\ref{A*eqn}) for every solution
$y$ of (\ref{b1}), where  $U'/U = -\nu_1$. The equation (\ref{b1}) has a pair of linearly independent 
canonical formal solutions with distinct exponential parts, and (\ref{A*eqn}) has an irregular singular point at infinity. 

If $\kappa$ is a non-zero exponential part for equation (\ref{b1}), then there exists 
a locally analytic solution $u_1$ of (\ref{b1}), with exponential part $\kappa$, which continues without zeros in some $\Omega(r_2)$ and 
annihilates a non-trivial locally analytic solution $y_1$ of $L[y]=0$, where $y_1$  is given by 
\begin{equation}
 \label{f1rep}
\frac{y_1'}{y_1} = d_0  \frac{u_1'}{u_1}  + d_1 .
\end{equation}
Moreover,  both $zu_1'(z^2)/u_1(z^2)$ and $zy_1'(z^2)/y_1(z^2)$ are rational at infinity.
\end{lem}
\textit{Proof.} The existence of the equation (\ref{A*eqn}) solved by $y U^{-1/2} $ for every solution
$y$ of (\ref{b1}) is a standard consequence of Abel's identity. Now the exponential parts $\kappa_1, \kappa_2$ for (\ref{b1}) are polynomials in $z^{1/2}$, by (\ref{hille*}), 
and their sum 
is a polynomial in $z$; thus  $\kappa_j (z) = Q_j(z) + z^{1/2} (-1)^j Q^*(z)$
with $Q^*$ and the $Q_j$ polynomials in $z$.

Suppose that $\kappa_1 = \kappa_2 =  \kappa_0$. Then $\kappa_0$ is a polynomial,
and must be non-constant since $g$ satisfies (\ref{b1})
and  $f/F = g^k$ has an essential singularity at infinity. But this implies the existence of a sector on which every solution
of (\ref{b1}), including $g$, tends to zero transcendentally fast as $z \to \infty$, which 
contradicts Lemma \ref{lem1}.
Thus $\kappa_1 \neq \kappa_2$, so that (\ref{A*eqn}) has an irregular singular point at infinity, and
at least one canonical formal solution of (\ref{b1}) has 
non-zero exponential part. 

Take a canonical formal solution $u_1$ of (\ref{b1}) with exponential part
$\kappa \neq 0$. Then $u_1$ is given by a formal expression as in (\ref{formalsol}), but free of logarithms, and $u_1'/u_1$ is a formal series in 
descending powers of $z^{1/2}$.
Since  $f$ has finite order and an unbounded sequence of poles, the function $g'/g$ is not
rational at infinity. Thus $g$ cannot solve a first order homogeneous 
linear differential equation with coefficients which are rational at infinity, and so  the division
algorithm for linear differential operators \cite[p.126]{Ince} shows that the operator $N$ of Lemma \ref{lemgeqn} 
must satisfy $N = N_0 \circ (D^2 +\nu_1 D + \nu_0)$, for some operator $N_0$ of order $1$ or $0$. 
Hence every solution of (\ref{b1}), including $u_1$,  solves $N[y]=0$. 
It follows from Lemma \ref{lemexppart} that $u_1$ annihilates some canonical formal solution $y_1$ of $L[y]=0$.
This gives, using (\ref{b1}), (\ref{5}) and (\ref{d2def}), 
$$
- y_1' u_1  = y_1 V[u_1] = y_1 ( (\beta - \alpha \nu_1) u_1'  +  (\gamma - \alpha \nu_0 ) u_1 ) =
y_1( - d_0 u_1'  - d_1 u_1) ,
$$
and hence (\ref{f1rep}). Thus  $y_1'/y_1$ is also a formal series in  $z^{1/2}$, and 
the hypotheses of Lemma \ref{lemb2} are satisfied with $p=y_1$ and $q=u_1$. Hence (\ref{b3}) holds 
with the $b_j$ rational at infinity and $b_k \not \equiv 0$ by Lemma \ref{d0lem}.  
But $L[y_1]=0$, and so $u_1'/u_1$ is algebraic at infinity, 
that is, $u_1'/u_1$ solves a polynomial equation with coefficients which are rational at infinity.
In particular,  the series for $u_1'/u_1$ converges for large $z$ in some sector, as does that for $y_1'/y_1$, by (\ref{f1rep}), and 
$u_1$ and $y_1$ are analytic local solutions of (\ref{b1}) and $L[y]=0$ respectively. 
Since the algebraic equation for $u_1'/u_1$
has only finitely many branches for its solutions, and each branch has no poles in some sector
$|z| > r_2 , \, | \arg z  | < 4 \pi  $, it follows that $u_1$ continues without zeros in $\Omega(r_2)$. 
This means that, as $z$ crosses a critical ray of (\ref{A*eqn}), the solution $u_1U^{-1/2}$ of (\ref{A*eqn})
must change from small to large or vice versa. Therefore continuing twice around a circle $|z| = r_3 > r_2$ 
brings $u_1U^{-1/2}$ back to a constant multiple of itself, and the same is true for $u_1$. 
Thus $zu_1'(z^2)/u_1(z^2)$ is  rational at infinity, and so is $zy_1'(z^2)/y_1(z^2)$ by (\ref{f1rep}).
\hfill$\Box$
\vspace{.1in}

Choose a critical ray $\arg z = \theta^*$ for the equation (\ref{A*eqn}) and a
sector $S^*$, symmetric about the critical ray, and with internal angle slightly less than 
$4 \pi /(2 + \deg_\infty A^*)$,
in which $f$ has an unbounded sequence of poles, these being zeros of $g$. 
In the sector $S^*$, equation (\ref{A*eqn}) has two linearly independent zero-free analytic
solutions, by (\ref{hille*}). Denote these by
$u^*=u U^{-1/2}$ and $v^* = vU^{-1/2}$ say, where $u$ and $v$ solve (\ref{b1}).
Here $u$ and $v$  have distinct exponential parts $\kappa_u $ and $\kappa_v$, each a polynomial in $z^{1/2}$, and  it may be assumed that 
$\kappa_u$ is non-constant and
\begin{equation}
 \label{uchoose}
\liminf_{z \to \infty, z \in S^*} \left| \frac{\kappa_v(z)}{\kappa_u(z)} \right| \leq 1 , \quad g = v-u ,
\end{equation}
since $u$ and $v$ may be interchanged and multiplied by constants.
Now Lemma \ref{lemb3*} shows that  there exist locally analytic 
solutions $u_1$ of (\ref{b1}) and  $y_1$ of $L[y]=0$ respectively, such that $u_1$ has exponential part $\kappa_u  $,
while (\ref{f1rep}) holds 
and both $zu_1'(z^2)/u_1(z^2)$ and $zy_1'(z^2)/y_1(z^2)$ are rational 
at infinity. Thus $u_1$ must be a constant multiple of $u$ and so, by (\ref{b0*}),  
\begin{equation}
 \label{f'frepa}
T_1 =  \frac{y_1'}{y_1} = d_0  \frac{u'}{u}  + d_1, \quad 
\frac{f'}{f} = d_0  \frac{g'}{g}  + d_1 =  d_0 \left( \frac{g'}{g} - \frac{u'}{u} \right) + \frac{y_1'}{y_1} .
\end{equation}
Poles $z$ of $f$ occur where $v(z)/u(z) = 1$, and have multiplicity 
equal to $-d_0(z)$, by (\ref{f'frepa}). 
Furthermore, by (\ref{hille*}),  $\zeta = (1/2\pi i) \log ( v^*/u^* ) = (1/2\pi i) \log ( v/u ) $ 
maps $S^*$ conformally onto a domain containing a right or left half-plane
$\pm {\rm Re} \, \zeta > M_0 > 0$. Since $d_0$ takes integer values at all points in $S^*$ 
where $\zeta$ is integer-valued, applying Lemma \ref{lemintegervalued} shows that there exists a polynomial $Q$ such that 
\begin{equation}
 \label{d0Q(M)}
d_0 = Q(T), \quad T = 2 \pi i \zeta = \log \left( \frac{v}{u} \right) .
\end{equation}
The second equation of (\ref{f'frepa}) can now be written  in the form
\begin{eqnarray}
 \frac{f'}{f} &=& Q(T) \left( \frac{v'-u'}{v-u} - \frac{u'}{u} \right) + \frac{y_1'}{y_1} =
Q(T) \frac{v'u-u'v}{(v-u)u}  + \frac{y_1'}{y_1} \nonumber \\
&=&
Q(T) \frac{v'/v-u'/u}{1-u/v }  + \frac{y_1'}{y_1} =
\frac{Q(T) T'}{1-e^{-T} }  + \frac{y_1'}{y_1} = \frac{Q(T) T'}{1-e^{-T} }  + T_1  , 
\label{f'/fQform}
\end{eqnarray}
which gives (\ref{caseIIrep}), and it suffices to consider two cases.

\subsection{Case I}

Suppose first that $Q$  is constant and one  exponential part for (\ref{b1}) is $0$. 
Then $d_0 = Q(T)$ is constant and 
$v$ has exponential part $0$ in $S^*$, because $u$ does not. A pole of $f$ of multiplicity
$m_0$ in $S^*$ gives $v/u=1 $ and  $g' = v'-u' = v'- (u'/u) v = T'v$, as well as 
(\ref{ac7}).
Since all poles of $f$ in $S^*$ with $|z|$ sufficiently large
have fixed multiplicity $-d_0$, it follows from Lemma \ref{lemintegervalued} and (\ref{ac7}) that $(T'v)^{-k}$,
which also has exponential part $0$ in $S^*$, must be constant, and so must $v'-(u'/u) v$. But then $W(u,v)/u$ is constant,
and so $\nu_1 = - u'/u$ in (\ref{b1}). Because $u$ solves (\ref{b1}), it must be the case that $\nu_0 = -(u'/u)' =  \nu_1'$. 
Now comparing (\ref{b1}) and
(\ref{newgeqn}) shows that, since $d_0$ is constant, $c_{k-2} - C_{k-2}$ must vanish, so that Proposition \ref{FHcase} may be applied,
and $f$ satisfies conclusion (i) of Theorem~\ref{thm1}.

\subsection{Case II}

Assume now that either both exponential parts for (\ref{b1}) are non-zero, or $Q$ is non-constant.

\begin{lem}
 \label{lemcaseII}
The solution $v$ continues without zeros in some $\Omega(r_2)$, and $z v'(z^2)/v(z^2)$ is rational at infinity. 
\end{lem}
\textit{Proof.}
Suppose first that both exponential parts for (\ref{b1}) are non-zero.
Then Lemma \ref{lemb3*} gives a solution $V_1$ of (\ref{b1}),
such that $V_1$ and $u$ are linearly independent and $V_1$ continues without zeros in some $\Omega(r_2)$,  
with $z V_1'(z^2)/V_1(z^2)$ rational at infinity. Since $u$ and $v$ are 
linearly independent and zero-free in $S^*$, the solution $V_1$ must be a constant multiple of $v$.

Now suppose that $v$ has exponential part $0$ in $S^*$: 
then  $Q$ is non-constant, and  (\ref{d0Q(M)}) implies
that  $T=\log (v/u)$ is algebraic at infinity.
Thus $v$ continues without zeros
in some $\Omega(r_2)$,
because $u$ does, and the same argument as  applied to $u$ in the proof of Lemma \ref{lemb3*} shows that $z v'(z^2)/v(z^2)$ is rational at infinity as asserted.
\hfill$\Box$
\vspace{.1in}

The functions $u'/u$, $v'/v$ and $T' = v'/v - u'/u$ are all defined for large $z \in S^*$ and given 
by convergent series in descending powers of $z^{1/2}$. Denote by $\widehat \psi $ the result of 
continuing a function element $\psi$ once counter-clockwise around a circle $|z|=r_3 > r_2$, starting in $S^*$. 
Since $u$ and $v$ both continue without zeros, 
there exists $\zeta_0 \in \C$ such that
\begin{equation}
 \label{continue}
 (a) \quad \widehat u = cu, \quad \widehat v = cv , \quad 
 \widehat T = T + \zeta_0 \quad \hbox{or}  \quad (b) \quad \widehat u = cv, 
\quad \widehat v = cu , \quad\widehat T = - T - \zeta_0 .
\end{equation}

\begin{lem}
\label{lemunotbig}
There exist $d_2 \in [0, 1/2]$ and  functions $E_0$, $E_1$ and $E_2 \not \equiv 0$, each rational at infinity, such that
\begin{equation}
 \label{uvrep1}
 \frac{u'(z)}{u(z)} = (d_2 -1) T'(z) + o( |T'(z)| ) , \quad 
\frac{v'(z)}{v(z)} = d_2  T'(z) + o( |T'(z)| ) 
\end{equation}
as $z \to \infty$ in $S^*$, while
$E_2 f'' + E_1 f' + E_0 f$ has no  zeros in some $\Omega(r_3)$. 

If subcase (a) applies in (\ref{continue}), then $T'$ is rational at infinity, with $T'(\infty) \neq 0$, while if  
subcase (b) applies
then $d_2 = 1/2$ and $H_0(z) = z^{1/2} T'(z)$ is rational at infinity, with $H_0(\infty) \neq 0$. 
\end{lem}
\textit{Proof.} 
Suppose first that subcase (a) applies in (\ref{continue}). 
Then $u'/u$, $v'/v$ and $T'$ are all rational at infinity, and 
so is $T_1$ in (\ref{f'frepa}). Thus
applying Lemma  \ref{lemb2} to $f$ and $g$ gives, in view of (\ref{b0*}), (\ref{f'/fQform}) and
Lemma \ref{d0lem2}, functions $E_0$, $E_1$ and $E_2$, each rational at infinity, such that $E_2 \not \equiv 0$ and 
$$
E_2 \frac{f''}{f} + E_1 \frac{f'}{f} + E_0 = \left( \frac{f'}{f} - T_1 \right)^2 
= \left(  \frac{Q(T) T'}{1-e^{-T} } \right)^2 .
$$
Hence $E_2 f'' + E_1 f' + E_0 f$ has no  zeros in some $\Omega(r_3)$. 

To prove the existence of $d_2$ in subcase (a), suppose first that 
$\deg_\infty (u'/u) > \deg_\infty T'$. Then, as $z \to \infty$, with $\arg z$ arbitrary, 
$$
\frac{v'(z)}{v(z)} = \frac{u'(z)}{u(z)} + T'(z) = (1+ o(1)) \frac{u'(z)}{u(z)} . 
$$
Since $u$ has non-zero exponential part, this gives a sector 
on which $u$ and $v$ both tend to zero transcendentally fast, 
and hence so does every solution of (\ref{b1}), including $g$, 
contradicting Lemma~\ref{lem1}.
Thus there exists $d_2 \in \C$ such that (\ref{uvrep1}) holds  as $z \to \infty$, with $\arg z$ arbitrary,
and $T'(\infty) \neq 0$, since $u$ has non-zero exponential part.
If $d_2 \not \in \R$, or if $d_2 \in \R \setminus [0, 1]$, then again 
there exists a sector 
on which $u$, $v$ and  $g$ all tend to zero transcendentally fast, contradicting Lemma \ref{lem1}. Finally, (\ref{uchoose}) gives
$d_2 \leq 1/2$.

Assume now that subcase (b) holds in (\ref{continue}). 
Because $f$ has an unbounded sequence of poles in $S^*$ and $y_1$ continues without zeros, 
(\ref{f'/fQform}) leads to
\begin{equation}
 \label{repcontd}
\widehat T' = - T', \quad 
 \frac{f'}{f}  =
\frac{-Q(T) T'}{1-e^{T+\zeta_0} }  + \widehat T_1 , \quad e^{\zeta_0} = 1, \quad  Q(T)T' = \widehat T_1 - T_1 =   T_2 - T_1 .
\end{equation}
Furthermore, $u'/u + v'/v = 2H_1$ and $u'v'/uv$ are rational at infinity,  and so are
$T_1+T_2$ and $T_1T_2$ by  continuation of the first equation of (\ref{f'frepa}).
On the other hand (\ref{repcontd}) implies that 
$T'(z) = 2 z^{1/2} H_2 (z)$, with $H_2$ rational at infinity.  
This yields 
\begin{eqnarray}
 \label{caseiib}
\frac{u'(z)}{u(z)} &=& H_1(z) - z^{1/2} H_2(z) , \quad \frac{v'(z)}{v(z)} = H_1(z) + z^{1/2} H_2(z) .
\end{eqnarray}
Since $u$ has non-zero exponential part, either $\mathrm{deg}_\infty \, H_1  \geq 0$
or $\mathrm{deg}_\infty \, H_2 \geq -1$. 
Moreover, $\mathrm{deg}_\infty \, H_2 \geq \mathrm{deg}_\infty \, H_1$ (and so $\mathrm{deg}_\infty \, H_2 \geq -1$) in (\ref{caseiib}); 
otherwise there again exists a sector on which 
every solution of (\ref{b1}), including $g$, tends to zero transcendentally fast, contradicting Lemma \ref{lem1}.
Thus (\ref{uvrep1})  holds with 
$d_2 =  1/2$.

Applying Lemma \ref{lemb2} to $f$ and $g$ now gives, in view of Lemma  \ref{d0lem2} and 
(\ref{b0*}),  (\ref{f'/fQform}) and (\ref{repcontd}), functions $E_0$, $E_1$ and $E_2$, 
each rational at infinity, such that $E_2 \not \equiv 0$ and 
$$
E_2 \frac{f''}{f} + E_1 \frac{f'}{f} + E_0 = \left( \frac{f'}{f} - T_1 \right) 
\left( \frac{f'}{f} - T_2 \right) =    - \frac{\left(Q(T) T'\right)^2}{ (1-e^T)(1-e^{-T})  }  ,
$$
and so $E_2 f'' + E_1 f' + E_0 f$ again has no zeros in some $\Omega(r_3)$. 
\hfill$\Box$
\vspace{.1in}

Recall that $\zeta(z) = T(z)/2 \pi i $ maps a subdomain  of $S^*$ conformally onto a right or left 
half-plane.
If $z_1 \in S^*$ and  $\zeta(z_1)   \in \Z$ then 
$e^{T(z_1)} = v(z_1)/u(z_1)=1$, while $f$ has a pole at $z_1$ 
of multiplicity $-d_0(z_1) = -Q(T(z_1))$, by (\ref{f'/fQform}), and  (\ref{ac7}) gives 
\begin{eqnarray}
 \label{g=0equ}
(T'(z_1)v(z_1))^{-k} &=& (v'(z_1)-u'(z_1))^{-k} =  g'(z_1)^{-k} =  Q_0(T(z_1))  , \nonumber\\
Q_0 &=& Q(Q-1) \ldots (Q-k+1). 
\end{eqnarray}
Lemma \ref{lemunotbig} makes it possible to write, on $S^*$, 
\begin{equation}
 \label{uvrep2}
T'(z)^k v(z)^k Q_0(T(z)) = e^{kd_2 T(z)} u_0(z) = e^{2 \pi i kd_2 \zeta (z)} u_0(z),  
\end{equation}
in which $ \log^+ |u_0(z)| = o( |T(z)| ) = o( |\zeta(z)| ) $ as $z \to \infty$ in $S^*$.
Thus 
(\ref{g=0equ}), (\ref{uvrep2}) and Lemma \ref{lemintegervalued} together imply that $k d_2 \in \Z$ and $u_0 \equiv 1$, so that
$v$ and $u$ have representations,
for some branch of $Q_0(T)^{1/k}$,   
\begin{equation}
\label{uvrep3}
 v =  \frac{e^{d_2 T}}
{ Q_0(T)^{1/k} T'}  ,  \quad   u = v e^{-T} = \frac{e^{(d_2-1) T}}
{ Q_0(T)^{1/k} T'} , \quad d_2 \in [0, 1/2], \quad kd_2 \in \Z, 
\end{equation}
and if $T'$ is not rational at infinity then $d_2 = 1/2$ and  $k$ is even. 
Now Abel's identity, (\ref{b1}), (\ref{newgeqn}),  (\ref{f'frepa}), (\ref{d0Q(M)}) and (\ref{uvrep3}) lead to
\begin{eqnarray*}
W_0 &=& W(u, v) = c e^{(2d_2-1) T} Q_0(T)^{-2/k} (T')^{-1} , \\
\nu_1  &=& \frac{d_0'+d_1 + a_{k-1}/k}{d_0 - (k-1)/2} =  
- \frac{W_0'}{W_0} = (1-2d_2)  T' + \frac{2 Q_0'(T)T'}{kQ_0(T)} + \frac{T''}{T'} ,  \\
\frac{y_1'}{y_1} 
&=& Q(T) \left( (d_2-1) T'   - \frac{Q_0'(T)T'}{k Q_0(T) } - \frac{T''}{T'} \right) +  \\
& & + \left( Q(T) - \frac{k-1}2 \right) \left((1-2d_2) T' + \frac{2 Q_0'(T)T'}{k Q_0(T) } + \frac{T''}{T'} \right) 
 - Q'(T)T' - \frac{a_{k-1}}k .
\end{eqnarray*}
Hence $T_3 = y_1'/y_1 + a_{k-1}/k$ is given by
\begin{equation}
 \label{fully1rep}
T_3 = \frac1k    \sum_{j=0}^{k-2} 
\left( \frac{j-k+1}{Q(T) - j} \right) Q'(T) T' -     
\left(   d_2 (Q(T) - k+1)  + \frac{k-1}2 \right) T'    - \left( \frac{k-1}2 \right)   
 \frac{T''}{T'}   .
\end{equation}
Thus 
$T_1 = y_1'/y_1 $ belongs to the field $\widetilde \Lambda$ generated by $d_0 = Q(T)$, $T'$, $a_{k-1}$ and their derivatives.
Since $L[y_1]=0$, a  standard change of variables gives
a linear differential operator $\widetilde L$ with coefficients $\widetilde c_j \in \widetilde \Lambda$ such that  
$$
L[y_1w] = y_1 \widetilde L[w], \quad \widetilde L =  \sum_{j=1}^{k} \widetilde c_j D^j , \quad \widetilde c_k = 1, \quad D = \frac{d}{dz} . 
$$
As $T_1$ is known,  the $ \widetilde c_j$  can be computed from the $ a_j $, and vice versa. Using (\ref{uchoose}) and (\ref{f'/fQform}), write 
\begin{eqnarray}
 \label{phidef1}
f &=& y_1 \phi, \quad \frac{1}{(v-u)^k} =
\frac1{g^k} = \frac{L[f]}f = \frac{\widetilde L[ \phi]}{\phi} =  \sum_{j=1}^{k} \widetilde c_j \, \frac{\phi^{(j)}}{\phi} , \nonumber \\
\frac{\phi'}{\phi} &=& \frac{S_1}{Y_1} , \quad S_1 = R_{1,0} = Q(T)T', \quad Y_1 = 1 - e^{-T}, \quad Y_1' =  T'(1-Y_1). 
\end{eqnarray}
There exist computable coefficients $ R_{j,\mu} \in \widetilde \Lambda $ such that, for $j \in \N$,
\begin{equation}
 \label{Rjmudef}
\frac{\phi^{(j)}}{\phi} = \frac{S_j}{Y_1^{j}} , \quad S_j =  \sum_{\mu=0}^{j-1} R_{j,\mu} Y_1^\mu , 
\quad R_{j,0} = Q(T)(Q(T)-1) \ldots ( Q(T)-j+1) (T')^j .
\end{equation}
The relations (\ref{Rjmudef}) hold by a straightforward induction argument, since the $S_j$ satisfy
$$
S_{j+1} = Y_1 S_j' - jY_1' S_j + S_1 S_j = Y_1 S_j' + j T'(Y_1-1) S_j + S_1 S_j,
\quad R_{j+1,0} = (Q(T)-j) T' R_{j,0} .
$$ 
Using (\ref{uvrep3}), (\ref{phidef1}) and (\ref{Rjmudef}) now delivers 
\begin{eqnarray*}
 Q_0(T) (T')^k (1-Y_1)^{kd_2} &=& Q_0(T) (T')^k e^{-kd_2T} = v^{-k} = \frac{Y_1^k}{(v-u)^k} 
= \sum_{j=1}^{k} \widetilde c_j S_j Y_1^{k-j}  \\
&=&   \sum_{j=1}^k \,  \sum_{\mu=0}^{j-1} \widetilde c_j R_{j,\mu} Y_1^{k-j+\mu } 
= \sum_{\mu=0}^{k-1} \sum_{j=\mu+1}^k \widetilde c_j R_{j,\mu} Y_1^{k-j+\mu } \\
&=&  \sum_{\mu=0}^{k-1} \,   \sum_{\nu=\mu}^{k-1} \widetilde c_{k+\mu-\nu}  R_{k+\mu-\nu,\mu} Y_1^{\nu }  
= \sum_{\nu=0}^{k-1} Y_1^{\nu }   \sum_{\mu=0}^{\nu} \widetilde c_{k+\mu-\nu}  R_{k+\mu-\nu,\mu}  ,
\end{eqnarray*}
in which $\nu = k-j+\mu$. Now $e^{-T} = u/v$ grows transcendentally fast on a subsector of $S^*$, 
whereas each element of the field $\widetilde \Lambda$ has form 
$\widetilde v (z) = v_1(z) + z^{1/2} v_2(z) $, with $v_1$ and $v_2$ rational at infinity.
Thus $e^{-T}$ is transcendental over 
$\widetilde \Lambda$ and so is $Y_1 = 1-e^{-T}$. 
Since (\ref{d0Q(M)}), (\ref{g=0equ}), (\ref{Rjmudef}) and Lemma~\ref{d0lem} together imply 
that $R_{j,0} \not \equiv 0$ for $j=1, \ldots, k$ and  that $\widetilde c_k R_{k,0} = Q_0(T) (T')^k $, 
comparing the coefficients of $Y_1^\nu$, starting 
from $\nu=1$, determines successively $\widetilde c_{k-1} , \ldots, \widetilde c_1$ and hence $\{ a_0, \ldots, a_{k-2} \}$.
\hfill$\Box$
\vspace{.1in}

The proof of Proposition~\ref{frankcase} is complete, but it is worth remarking that  (\ref{deta1}),  (\ref{exL1})  and (\ref{exL2}) 
show that $d_2 =0$ and $d_2 = 1/3$ are both possible when $k=3$. 
Furthermore, Propositions \ref{FHcase}
and ~\ref{frankcase} each give a solution $y_0$ of $L[y]=0$ on $S^*$, of the form (\ref{analyticsol}), whose exponential part is a non-constant polynomial 
in $z^{1/2}$, and 
in $z$ if $k$ is odd. If  $d_2 = 0$ or $d_2=1/2$ in (\ref{fully1rep}), or if $f$ is given by (\ref{Brurep}), 
then $y_0 = y_1$,  by (\ref{d0Q(M)}), Lemma \ref{lemunotbig} and the fact that $H''/H'$ does not vanish at infinity in  (\ref{Brurep}).
On the other hand, if $0 < d_2 < 1/2$ then $u$ and $v$ both have
non-constant exponential part, by Lemma \ref{lemunotbig} and (\ref{uvrep3}), and Lemma \ref{lem1} 
gives a non-trivial solution $y_2$ of $L[y]=0$  annihilated by $v$; thus  
$y_2'/y_2 - y_1'/y_1 = d_0 (v'/v - u'/u) = Q(T)T'$, by (\ref{d2def}), (\ref{f'frepa}) and
(\ref{d0Q(M)}),
and $y_0 \in \{ y_1, y_2 \}$. 
It follows that $y_0'(z)/y_0(z) = (1+o(1))c z^{\lambda_0} $ and
$y_0^{(j)}(z)/y_0^{(k)}(z) = (1+o(1))c z^{(j-k)\lambda_0} $ as $z \to \infty$ in a subsector of $S^*$, for $j = 0, \ldots, k-1 $,
where $\lambda_0 \geq -1/2$,
and $\lambda_0 \geq 0$ if $k$ is odd. This implies that at least one $a_j$ has $a_j(z) \neq  o( |z|^{(k-j) \lambda_0} ) $
as $z \to \infty$, which is sharp by \cite[(1.8)]{Latsuji}.

\section{A change of variables}\label{changevariables}

In order to prove Theorem \ref{thm1} it now suffices, in view of Proposition \ref{FHcase}, to 
show that the hypotheses of Proposition \ref{frankcase} are satisfied when $L \neq M$. 
Since the value of $E$ is immaterial in Proposition \ref{frankcase}, 
a change of variables $z \to z^n$ may now be employed 
to ensure that, in the terminology of Section \ref{wasowthm}, the integer $p$ is $1$, so
that the exponential parts  and associated asymptotic or formal series involve only integer powers. 
Indeed, let $k \geq 3$ and $n \geq 2$ be integers and let $f$, $F$ and $\mathfrak{f}$ satisfy 
\begin{equation}
 \label{cf1}
F(z) = L[f](z) = f^{(k)}(z) + a_{k-1}(z) f^{(k-1)}(z) + \ldots + a_0(z) f(z) , \quad \mathfrak{f}(z) = f(z^n) , 
\end{equation}
where the $a_j$ are rational at infinity with $a_{k-1}(\infty ) = 0$. Take linearly independent locally analytic solutions
$f_1, \ldots, f_k$ of $L[y]=0$.

\begin{lem}
 \label{lemcf1} For each integer $m \geq 1$ there exist  rational functions 
$c_{p,m}(z)$, depending only on $m$ and $n$, such that 
\begin{equation}
 \label{cf3}
f^{(m)}(z^n) = \sum_{p=1}^m c_{p,m} (z) \mathfrak{f}^{(p)}(z), \quad  c_{m,m}(z) = (nz^{n-1})^{-m} .
\end{equation}
Moreover, if $m \geq 2$ then $c_{m-1,m} (z) /c_{m,m}(z) \to 0$ as $z \to \infty$. 
\end{lem}
\textit{Proof.}
Clearly $\mathfrak{f}(z) = f(z^n)$  gives, as $z \to \infty$, 
$$
f'(z^n) = (nz^{n-1})^{-1} \mathfrak{f}'(z), \quad f''(z^n) = (nz^{n-1})^{-2} \mathfrak{f}''(z) 
+ c_{1,2}(z) \mathfrak{f}'(z),
\quad c_{1,2}(z) = O( |z|^{1-2n}) .
$$
Next, if the assertions of the lemma hold for some $m \geq 2$ then, as $z \to \infty$, 
\begin{eqnarray*}
 nz^{n-1} f^{(m+1)}(z^n) &=& c_{m,m}(z) \mathfrak{f}^{(m+1)}(z) + \mathfrak{f}^{(m)}(z) ( c_{m,m}'(z) + c_{m-1,m}(z)) + \ldots \\
&=&
c_{m,m}(z) \left[ \mathfrak{f}^{(m+1)}(z) + \mathfrak{f}^{(m)}(z)  O( 1/|z| ) + \ldots \right] .
\end{eqnarray*}
\hfill$\Box$
\vspace{.1in}

Now (\ref{cf1}) and (\ref{cf3}) yield,
as $z \to \infty$, 
\begin{eqnarray*}
F(z^n)  &=&  (nz^{n-1})^{-k} \left[ \mathfrak{f}^{(k)}(z) + \mathfrak{f}^{(k-1)}(z) O( 1/|z| ) + \ldots \right] + \\
&& + O( |z|^{-n} ) (nz^{n-1})^{1-k} \left[ \mathfrak{f}^{(k-1)}(z)  + \ldots \right] + \ldots \\
&=& (nz^{n-1})^{-k} \left[ \mathfrak{f}^{(k)}(z) + \mathfrak{f}^{(k-1)}(z) O( 1/|z| ) + \ldots \right] .
\end{eqnarray*}
Hence there exist  functions $\mathfrak{a}_j(z)$, all rational at infinity
and with $\mathfrak{a}_{k-1}(\infty) = 0$, such that 
\begin{equation*}
 \label{cf4}
\mathfrak{f}^{(k)}(z) + \mathfrak{a}_{k-1}(z) \mathfrak{f}^{(k-1)}(z) + \ldots + \mathfrak{a}_0(z) \mathfrak{f}(z) = \mathfrak{F}(z) = (nz^{n-1})^k F(z^n) .
\end{equation*}
The new operator is $\mathfrak{L} = D^k + \mathfrak{a}_{k-1} D^{k-1} + \ldots $, 
and $y(z^n)$ solves $\mathfrak{L}[y] =0$ for every locally analytic or formal solution $y$ of $L[y]=0$, as does each $\mathfrak{f}_j (z) = f_j(z^n)$.

If $f$ is as in the hypotheses of Theorem \ref{thm1} then running Frank's method as in Section \ref{frankmethod} for $\mathfrak{f}$ and
$\mathfrak{F}$ gives rise to auxiliary functions $\mathfrak{g}$,  $\mathfrak{h} = - (\mathfrak{f}'/\mathfrak{f})\mathfrak{g}$
and $\mathfrak{w}_j$,  which 
satisfy, using (\ref{2}), 
\begin{eqnarray*}
\mathfrak{g}(z)^k &=& \frac{\mathfrak{f}(z)}{\mathfrak{F}(z)} = \frac{f(z^n)}{(nz^{n-1})^k F(z^n) } 
= \frac{g(z^n)^k}{(nz^{n-1})^k  } 
, \quad 
 \mathfrak{g}(z) = \frac{g(z^n)}{n z^{n-1} } , \\
\mathfrak{h}(z)&=& - \frac{\mathfrak{f}'(z) }{\mathfrak{f}(z)} \mathfrak{g}(z) = -
n z^{n-1} \frac{f'(z^n)}{f(z^n)} \, \frac{g(z^n)}{n z^{n-1} } = h(z^n), \\
\mathfrak{w}_j(z) &=&
\mathfrak{f}_j' (z) \mathfrak{g}(z) + \mathfrak{f}_j(z) \mathfrak{h}(z)
= f_j'(z^n) g (z^n) + f_j(z^n) h(z^n) = w_j(z^n) . 
\end{eqnarray*}
Thus the $\mathfrak{w}_j$  solve the   equation
$\mathfrak{M}[y] = y^{(k)} + \ldots =0$ which is obtained from $M[y]=0$ in the same way as $\mathfrak{L}[y] =0$ arose from $L[y]=0$.
In $\mathfrak{M}[y]$ the coefficient of $y^{(k-1)}$ is $\mathfrak{a}_{k-1}$,
since $a_{k-1} = A_{k-1}$.
It is important to note that $\mathfrak{L} = \mathfrak{M}$ if and only if $L=M$.

Therefore $n$ may be chosen so that in the canonical formal solutions
(\ref{formalsol}) for the equations $\mathfrak{L}[y] =0$ and  $\mathfrak{M}[y] =0$ the integer $p$ is $1$.
Moreover, $\mathfrak{L}[y] =0$  has linearly independent canonical formal solutions  $\mathfrak{h}_1, \mathfrak{h}_2$ whose formal Wronskian
$W( \mathfrak{h}_1, \mathfrak{h}_2)$ is free of logarithms. This is clear if there are solutions of $\mathfrak{L}[y] =0$
with distinct exponential parts.
On the other hand if all exponential parts for $\mathfrak{L}[y] =0$ are the same then they are all $0$, since $\mathfrak{a}_{k-1}(\infty) = 0$,
and there exists a solution  $\mathfrak{h}_1 (z) = z^{e_1} R_1(z)\not \equiv 0$, with $e_1 \in \C$ and $R_1 $ rational at infinity. The standard
reduction of order method then gives an equation which is solved by $(y/\mathfrak{h}_1)'= W(\mathfrak{h}_1, y)\mathfrak{h}_1^{-2}$, 
for every solution $y$ of $\mathfrak{L}[y] =0$, and which has a canonical formal solution  free of logarithms.

When $L\neq M$, and hence $\mathfrak{L} \neq \mathfrak{M}$, Lemma \ref{lemgeqn} applied to
$\mathfrak{L} $ and $ \mathfrak{M}$ gives an 
equation $\mathfrak{N}[y]=0$, of order $2$ or $3$, which is solved by $\mathfrak{g}$, as well as a counterpart
$\mathfrak{V}$ for the operator $V$. 
Choosing $\mathfrak{h}_1, \mathfrak{h}_2$ as in the previous paragraph and any canonical formal solution $\mathfrak{G}$ of $\mathfrak{N}[y]=0$ then makes 
$\mathfrak{h}_j' \mathfrak{G} + \mathfrak{h}_j \mathfrak{V} [\mathfrak{G}] $ a solution of $\mathfrak{M}[y] =0$. Solving for $\mathfrak{G}$ by Cramer's rule shows that 
in the canonical formal solutions
(\ref{formalsol}) for  $\mathfrak{N}[y] =0$ it may also be assumed that $p=1$.

\section{The main step}\label{brug1}

\begin{prop}
 \label{mainprop}
 Assume that $f$ and $F$ are as in 
the hypotheses of Theorem \ref{thm1}, and that  
$f$ has finite order and an unbounded sequence of poles. Then there exist $E \in \N$ and a 
function $R$ which 
is rational
at infinity such that all  poles $z$ of $\mathfrak{f}(z) = f\left(z^E\right)$ with $|z| $ sufficiently large have multiplicity $R(z)$. 
\end{prop}

Once Proposition \ref{mainprop} is proved, Theorem \ref{thm1} is established as follows. 
Under the hypotheses of Theorem \ref{thm1}, 
the first possibility is that $L=M$ and $f$ is determined by 
Proposition \ref{FHcase}. If this is not the case then Lemma \ref{lemgeqn} shows that 
$f$ has finite order and an unbounded sequence of poles. In view of Proposition \ref{mainprop},
Proposition \ref{frankcase} may be applied, and $f$ is thereby determined. 
\hfill$\Box$
\vspace{.1in}

Assume for  the remainder
of the paper that the assumptions of Proposition \ref{mainprop} are satisfied. 

\begin{lem}
\label{lemmainstep}
 The following additional assumptions may all be made.\\
(A)  In Section \ref{frankmethod}, the operators $L$ and $M$ are not the same.\\
(B) In  the canonical formal solutions (\ref{formalsol}) for all of the equations
$L[y]=0$, $M[y]=0$ and $N[y]=0$, the integer $p$ is $1$.\\
(C) The function $g$ solves no second order homogeneous linear differential equation
with coefficients which are rational at infinity. Moreover, the operator $N$  has order $3$ and may be written in the form
\begin{equation}
 \label{Nform}
 N = D^3 + B_2 D^2 + B_1 D + B_0,
\end{equation}
with the $B_j$ rational at infinity,
while $\alpha \not \equiv 0$ in (\ref{5}) and $D_{k-2} = C_{k-2} - c_{k-2} \not \equiv 0$ in Section \ref{frankmethod}.
\end{lem}
\textit{Proof.}
Assumption (A) is legitimate because of  Proposition \ref{FHcase}, while (B)
is justified by taking $\mathfrak{f}(z) = f(z^{m_1})$ in place of $f$, 
for some $m_1 \in \N$, as in Section \ref{changevariables}. 
Next, the first three assumptions of (C) are valid since otherwise (\ref{5}) shows that $f$ and $g$ satisfy an equation (\ref{b0*}) with $d_0$ and $d_1$ rational at 
infinity, in which case the conclusion of
Proposition \ref{mainprop} follows from a comparison of residues. 
The last assumption of (C) is justified by Proposition \ref{FHcase}. 
\hfill$\Box$
\vspace{.1in}

\begin{lem}
 \label{propf'not0}
Assume that there exists  a  function
$a^*$ which is rational at infinity, with the property that $- f'/f  + d g'/g - a^*$ has no zeros in some $\Omega(r_2)$, where 
$d \in \{ 0, \ldots, k-1\}$ is a constant. Then $d$ satisfies $d \neq (k-1)/2$.

Assume further that $a^*(\infty) \neq 0$. Then 
$g$ is given by 
\begin{equation}
 \label{newestgrep}
P' g = \beta_1 e^{\omega_1 P} + \beta_2 e^{\omega_2 P} + \beta_3 e^{\omega_3 P} , \quad \beta_j, \omega_j \in \C^* , \quad 
1 = \omega_1 \neq \omega_2 \neq \omega_3 \neq 1 ,
\end{equation}
in which 
$P'$ is rational at infinity, with $P'(\infty) \neq 0$.
If, in addition, $d=0$ or $d=k-1$ then $f$ satisfies the conclusion of Proposition \ref{mainprop}.
\end{lem}
\textit{Proof.} 
As in  (\ref{gf1}), write $p'/p = f'/f + a_{k-1}/k$ and $q = -(p'/p)g$. Then $g$ and $q$ solve the equations (\ref{9818}) to
(\ref{D*eqnQ}). Moreover, $a = a^* - a_{k-1}/k$ is rational at infinity and $ -p'/p + dg'/g - a  $ has no 
zeros in $\Omega(r_2)$. 
Since poles of $p'/p$ have negative residues and 
are  simple zeros of $g$, while $f'/f$ and $g^k$ have finite order,
it is possible to write 
\begin{equation}
 \label{ePdef}
q + dg' -ag = 
 g \left( - \frac{p'}{p} + \frac{dg'}{g} - a \right) = e^P,
\end{equation}
with $P'$ rational at infinity.  
Then (\ref{9821X}) and (\ref{ePdef})  yield 
\begin{equation}
 \label{alt2a1}
P'e^P = x  g'' - a g' - \left( \frac{D_{k-2}}k + a' \right)  g, \quad  x =  d -  \frac{k-1}2 ,
\end{equation}
and by Lemma \ref{lemmainstep}(C) it may be assumed that $P' \not \equiv 0 $. 
Differentiation of this equation leads to 
\begin{eqnarray}
 \label{alt2a1a}
0 &=& x g''' + g'' \left( -x  \left( \frac{P''}{P'} +P' \right) - a \right) 
 +  g' \left( a \left( \frac{P''}{P'} +P' \right) - \frac{D_{k-2}}k - 2 a' \right) \nonumber \\
& &  + g \left( \left( \frac{P''}{P'} +P' \right) \left( \frac{D_{k-2}}k + a' \right)  - \frac{D_{k-2}'}k - a'' \right) ,
\end{eqnarray}
and so $x \neq 0$ and $d \neq (k-1)/2$, as asserted, again by Lemma \ref{lemmainstep}(C).

Now assume that $a^*(\infty) \neq 0$, which implies that  $a(\infty) \neq 0$.
The following  is an extension of a method from  \cite{GL}.  Since $C_{k-2} = D_{k-2} + c_{k-2} $, formula 
(\ref{9821Y}) 
becomes, in view of (\ref{ePdef}),
\begin{eqnarray}
 \label{alt2a2} 
D_{k-2}  e^P &=& \frac{k(k^2-1)}{12} g''' + g' \left( \left( x+1 \right)  D_{k-2} + 2 c_{k-2} \right)  \nonumber \\
& & +  g \left( \frac{k-1}2 D_{k-2}' + c_{k-2}' - D_{k-3} - a D_{k-2} \right) ,
\end{eqnarray}
and (\ref{D*eqnQ}) may be written as 
$$
 \left( \frac{2 D_{k-3}}{k-2} - D_{k-2}' \right) e^P =
\frac{(k+2) D_{k-2}}3 \,  g'' + d_5 g' + d_6 g,
$$
with $d_5, d_6$  rational at infinity. 
Comparing the last equation with (\ref{alt2a1})  delivers 
\begin{equation}
 \label{alt2a3}
x  \left( \frac{2 D_{k-3}}{k-2} - D_{k-2}' \right) =  \frac{(k+2) D_{k-2} P'}3 ,
\end{equation}
again using  Lemma \ref{lemmainstep}(C).
Combining (\ref{alt2a1}) with (\ref{alt2a2}) and  (\ref{alt2a3})  leads to
\begin{eqnarray}
\label{alt2a4}
 0 &=& \frac{k(k^2-1)}{12} g''' - g'' \left( \frac{x D_{k-2}}{P'} \right) 
+ g' \left( (x+1) D_{k-2} + 2 c_{k-2} + \frac{aD_{k-2}}{P'} \right) \nonumber \\
& & + g \left( \frac{D_{k-2}'}2 + c_{k-2}' - \frac{(k^2-4) D_{k-2}P'}{6x }  - a D_{k-2} + \frac{D_{k-2}}{P'} \left( \frac{D_{k-2}}k + a' \right) \right) .
\end{eqnarray}
Lemma \ref{lemmainstep}(C) implies that (\ref{alt2a4}) must be (\ref{alt2a1a}) multiplied by $k(k^2-1)/12x $. Comparing the coefficients of $g''$ yields 
\begin{equation}
 \label{alt2a5}
\frac{P''}{P'} +P' = - \frac{a}{x} + \frac{12x  D_{k-2}}{k(k^2-1) P'} .
\end{equation}
Next, matching the coefficients of $g'$ and using (\ref{alt2a5}) results in 
\begin{eqnarray}
 \label{alt2a6}
c_{k-2} 
&=& - \frac{(12x^2+12x+k^2-1)D_{k-2}}{24x}  - \frac{k(k^2-1)a'}{12x} - \frac{k(k^2-1)a^2}{24 x^2 } .
\end{eqnarray}
Examining the coefficients of $g$ in (\ref{alt2a1a}) and (\ref{alt2a4}) in the light of (\ref{alt2a5}) and (\ref{alt2a6}) leads to
$$
a D_{k-2} \left( \frac{k^2-1-12x^2}{12x^2} \right)  = D_{k-2}' \left( \frac{12x^2+1-k^2}{24x} \right) + \frac{(k^2-4)D_{k-2} P'}{6x} .
$$
Because $k \geq 3$ and $D_{k-2} P' \not \equiv 0$, this forces $k^2-1-12x^2 \neq 0$ and 
\begin{equation}
 \label{alt2a7}
a = - \, \frac{x D_{k-2}'}{2 D_{k-2}} + \frac{2x (k^2-4)P'}{k^2-1 -12x^2}  .
\end{equation}
Therefore $P'(\infty) \neq 0$, since $a(\infty) \neq 0$, and using (\ref{alt2a7}) to eliminate $a$ from (\ref{alt2a5}) delivers 
$$
\frac{P''}{P'} + \frac{(3k^2-12x^2-9) P' }{k^2-1-12x^2} - \frac{ D_{k-2}'}{2  D_{k-2}} = \frac{12x D_{k-2}}{k(k^2-1) P' } .
$$
Setting $Z = 1/D_{k-2}$ yields in turn a linear differential equation of form
\begin{equation}
 \label{alt2a8}
\left(  
\frac{2P''}{P'} + \eta_1 P' \right) Z + Z' =  \frac{24x   }{k(k^2-1) P' }   , \quad \eta_1 = \frac{6(k^2-4x^2-3)  }{k^2-1-12x^2}.
\end{equation}
If  $k^2  -4x^2 - 3 = 0$, then   $(k+2x)(k-2x) = 3$; because $k, 2x \in \Z$, this forces  $2k = \pm 4$, a contradiction. 
Assume henceforth that $k^2-4x^2-3 \neq 0$:
then the integrating factor for (\ref{alt2a8}) is $(P')^2 e^{\eta_1 P}$, with $\eta_1 \neq 0$, and the general solution to (\ref{alt2a8}) is 
$$
Z = (P')^{-2} \left( \eta_2 + d_7 e^{- \eta_1 P} \right), \quad d_7 \in \C, \quad \eta_2 = \frac{24x}{k(k^2-1) \eta_1 }. 
$$
Since 
$P'$ and $Z$ are rational at infinity with $P'(\infty) \neq 0$, this yields $Z = \eta_2 (P')^{-2} $ and 
\begin{equation}
 \label{eta3def}
\frac{-D_{k-2}}{k}  =   \frac{-(P')^2}{k \eta_2} = \eta_3 (P')^2, \quad \eta_3 = \frac{- (k^2-1) \eta_1 }{24x}
= \frac{(4x^2 +3 -  k^2 ) (k^2-1) }{4x(k^2-1-12x^2)} \neq 0 ,
\end{equation}
as well as
\begin{equation}
 \label{eta4def}
a = -x \left( \frac{ P''}{P'} + \eta_4 P' \right), \quad 
\eta_4 =  - \, \frac{2 (k^2-4)}{k^2-1 -12x^2} \neq 0,
\end{equation}
using (\ref{alt2a7}). Combining   (\ref{alt2a1}), (\ref{eta3def}) and (\ref{eta4def})  shows that $g$ solves the equation
\begin{equation}
 \label{alt2a10}
P'e^P = x  \left(   y' + \left( \frac{P''}{P'} + \eta_4 P'  \right) y \right)' + \eta_3 (P')^2  y   .
\end{equation}
Now write $\zeta = P(z)$ and $Y_0(\zeta) = y(z) P'(z)$ so that 
$$
y' + \left( \frac{P''}{P'} + \eta_4 P'  \right) y = \frac{dY_0}{d \zeta } + \eta_4 Y_0, \quad 
\left(   y' + \left( \frac{P''}{P'} + \eta_4 P'  \right) y \right)' = P' \left( \frac{d^2Y_0}{d \zeta^2 } + \eta_4 \frac{dY_0}{d \zeta } \right) .
$$
Thus (\ref{alt2a10}) becomes 
\begin{equation}
 \label{alt2a10a}
e^{\zeta} = x Y_0''(\zeta) + x \eta_4 Y_0'(\zeta) + \eta_3 Y_0( \zeta ).
\end{equation}
The auxiliary equation for the complementary function of (\ref{alt2a10a}) is then 
\begin{equation}
 \label{auxeq}
x \lambda^2 + x \eta_4 \lambda + \eta_3 = 0 , \quad x \eta_3 \eta_4 \in \C^*.
\end{equation}

Suppose that $1$ is a double root of (\ref{auxeq}). Then $g$ has a representation 
$g = (P')^{-1} e^P Q_2(P)$, for some polynomial $Q_2 \not \equiv 0$ of degree at most two. Since 
$P'(\infty) \neq 0$, there cannot exist a sector on which $g$ has an unbounded sequence of zeros, contradicting the assumption that 
$f$ has an unbounded sequence of poles.

Now suppose that $1$ is a simple root of (\ref{auxeq}), or that (\ref{auxeq}) has a repeated root. Then the fact that $\eta_3 \neq 0$ gives
$$
P'g = ( \beta_1 +  \beta_2 P ) e^{\omega_1 P} +  \beta_3 e^{\omega_3 P} , \quad \beta_j, \omega_j \in \C, \quad 0 \neq \omega_1 \neq \omega_3 \neq 0.  
$$
Here $ \beta_2 \neq 0$  by Lemma \ref{lemmainstep}(C), since otherwise $g$ satisfies a second order linear differential
equation,
and $ \beta_3 \neq 0$ by the assumption that $f$ has an unbounded sequence of poles. 
Denote by $\widehat \psi $ the result of analytically continuing a function element $\psi$ once around a given circle $|z| = r_3 > r_2$. Then 
there exists $\zeta_0 \in \C$ such that $\widehat P = P + \zeta_0$ and 
$$
(\beta_1 + \beta_2  \zeta_0 +  \beta_2 P ) e^{\omega_1 P +  \omega_1 \zeta_0} +  \beta_3 e^{\omega_3 P +  \omega_3 \zeta_0} =
 P' \, \widehat g =   P' \omega g =   \omega \left( ( \beta_1 +  \beta_2 P ) e^{\omega_1 P} +  \beta_3 e^{\omega_3 P} \right)  ,
$$
where $\omega^k =1 $. 
Because $\beta_3 P'(\infty) \neq 0$, examining the coefficients of $e^{\omega_3 P}$ and $ e^{\omega_1 P}$ leads to
$$e^{ \omega_3 \zeta_0} = \omega, \quad 
(\beta_1 +  \beta_2 \zeta_0 +  \beta_2 P ) e^{  \omega_1 \zeta_0 } = \omega  ( \beta_1 +  \beta_2 P ) .
$$
Differentiating the last relation then shows that  $ e^{\omega_1 \zeta_0 } = \omega $, since $\beta_2 \neq 0$,  and 
$$
\omega \beta_1 =  (\beta_1 +  \beta_2 \zeta_0 )  e^{ \omega_1 \zeta_0 } = \omega (\beta_1 +  \beta_2 \zeta_0 ) ,
$$
so that $\zeta_0 = 0$ and $P$ is rational at infinity, which forces $g$ to solve  a second order 
equation, contradicting Lemma \ref{lemmainstep}(C).
Thus (\ref{newestgrep}) holds, with the $\omega_j \in \C^*$ pairwise distinct, since $\eta_3 \neq 0$, and $\omega_1=1$,
and none of the $\beta_j$ can vanish, again by  Lemma \ref{lemmainstep}(C).
The proof of (\ref{newestgrep}) is now complete.

Next, suppose that  $d=0$ or $d=k-1$, so that $x = \pm (k-1)/2$. 
Now (\ref{eta3def}) and (\ref{eta4def}) imply that (\ref{alt2a10a}) takes the form
$$
e^{\zeta} = x Y_0''(\zeta) + x \left( \frac{k+2}{k-1} \right) Y_0'(\zeta) + x \left( \frac{k+1}{(k-1)^2} \right) Y_0( \zeta ).
$$
The auxiliary equation for the complementary function has  roots $\lambda_j = 1 - jk/(k-1)$, for $j=1, 2$, 
and (\ref{newestgrep}) becomes, in view of Lemma \ref{lemmainstep}(C),  
\begin{equation}
 \label{alt2a11}
g = \frac{e^P}{P'} \left( e_1 + e_2 e^{\eta P} + e_3 e^{2 \eta P} \right), \quad \eta = - \, \frac{k}{k-1}, \quad e_j \in \C^*.  
\end{equation}
Thus (\ref{gf1}), (\ref{ePdef}), (\ref{alt2a11}) and partial fractions deliver 
\begin{eqnarray}
 \label{alt2a12}
g &=& \frac{e_4 e^P}{P'} \left( e^{\eta P} - e_5 \right) \left( e^{\eta P} - e_6 \right) ,  \nonumber \\
\quad \frac{f'}{f} + \frac{a_{k-1}}k + a &=&  
\frac{p'}{p} + a 
 =  e_7 P'  \left( \frac{1}{e^{\eta P} - e_5 }  - \frac{1}{e^{\eta P} - e_6 } \right) + \frac{dg'}{g} , \quad e_j \in \C.
\end{eqnarray}
Again the  $e_j $ are all non-zero, and $e_5 \neq e_6$ since $g$ cannot have multiple zeros. 
If $r_4$ is large and some continuation of $e^{\eta P}$ takes the value $e_5$ at some $z_0 \in \Omega(r_4)$ 
then $z_0$ is a zero of $g$, and so is
a pole of $f$ of multiplicity $m_1 $ satisfying $- m_1 = d +  e_7 / \eta e_5 $, so that  (\ref{ac7}),  (\ref{alt2a11}) and (\ref{alt2a12}) imply that
at the point $z_0$ the following equations are satisfied: 
\begin{eqnarray*}
 e_5^{k-1} &=& e^{(k-1) \eta P} = e^{-kP} ; \\
\frac{(-1)^k}{ m_1 (m_1+1) \ldots (m_1+ k - 1) } 
&=& 
 (g')^k = \left(e_4 e^P \eta e_5 (e_5 - e_6 ) \right)^k \\
&=&  e_4^k  e_5^{1-k} \eta^k e_5^k (e_5 - e_6)^k   =  e_4^k e_5 \eta^k  (e_5 - e_6)^k  .
\end{eqnarray*}
Similarly, 
all zeros of continuations of  $e^{\eta P} - e_6$ to $\Omega(r_4)$ are poles of $f$ of multiplicity $m_2$, 
where $-m_2 = d - e_7/\eta e_6$, and 
$$
\frac{ m_2 (m_2+1) \ldots (m_2+ k - 1) }{ m_1 (m_1+1) \ldots (m_1+ k - 1) } = (-1)^k \, \frac{e_5}{e_6} = (-1)^{k-1} \frac{m_2+d}{m_1+d} . 
$$
But $d=0$ or $d=k-1$,  so that  $m_1 = m_2$ (and $k$ is odd). 
\hfill$\Box$
\vspace{.1in}

\section{The exponential parts for the equation $N[y]=0$}\label{expparts}

Let $p$ and $q$ be polynomials in $z$, and let $\theta \in \R$. 
Write $p \prec q$ (respectively, $p \preceq q$, $p \simeq q$) to indicate that 
${\rm Re} \, p( r e^{i \theta} ) < {\rm Re} \, q( r e^{i \theta} )$ 
(respectively ${\rm Re} \, p( r e^{i \theta} ) \leq  {\rm Re} \, q( r e^{i \theta} )$,
${\rm Re} \, p( r e^{i \theta} ) =  {\rm Re} \, q( r e^{i \theta} )$) as $r \to + \infty$. 
Since each $P_\theta (r) = {\rm Re} \, p( r e^{i \theta} )$ is a polynomial in $r$, every $\theta \in \R$
has $p \prec 0$ or $p \simeq 0$ or $0 \prec p$, and 
if $p $ is not constant then all but finitely many $\theta \in [0, 2 \pi ]$ 
have either $p \prec 0$ or $0 \prec p$.

Suppose that $N[y]=0$ has linearly independent canonical formal
solutions  with exponential parts
$\kappa_1, \kappa_2, \kappa_3$. The $\kappa_j$  are polynomials in $z$ by Lemma \ref{lemmainstep}, and it will be assumed 
as before that $\kappa_j (0) = 0$ for all $j$, from which it follows that if
$\kappa_j - \kappa_{j'} $ is constant then $\kappa_j -  \kappa_{j'} \equiv 0$.

\begin{lem}
 \label{lemkappanotsame}
The $\kappa_j$ are not all the same polynomial, and there does not exist $\theta \in \R$ with
$\kappa_j \prec 0$ on $\arg z = \theta$ for $j=1, 2, 3$. 
\end{lem}
\textit{Proof.} 
The first assertion is proved in Lemma \ref{lem1}, 
and the second holds 
because otherwise $g^k = f/F$ tends to zero transcendentally fast on a sector,  contradicting Lemma \ref{lem1}.
\hfill$\Box$
\vspace{.1in}

Lemma \ref{lemkappanotsame} does not exclude two of the $\kappa_j$ being the same polynomial, possibly identically zero, and this case will be
dealt with in Sections \ref{repeatnontriv} and \ref{repeattriv}. When there is no repetition among the $\kappa_j$,
the next lemma shows that there are two subcases to handle. 

\begin{lem}
 \label{lemkappaalldifferent}
Suppose that the $\kappa_j$ are pairwise distinct. 
Then it is possible to label the $\kappa_j$ and choose a ray
$\arg z = \theta \in \R$ such that
\begin{equation}
 \label{different1}
(A) \quad 
\kappa_1 \prec \kappa_2 \prec 0 \prec \kappa_3  \quad \hbox{or} \quad (B) \quad \kappa_1 \prec \kappa_2 \prec 0 , \quad \kappa_3 \equiv 0 
\end{equation}
or 
\begin{equation}
 \label{different2}
(C) \quad 
\kappa_1 \prec  0 \prec \kappa_3 , \quad \kappa_2 \equiv 0. 
\end{equation}
\end{lem}
\textit{Proof.} 
If one of the $\kappa_j$ is identically zero label the other two as $\kappa_a$ and $\kappa_b$, and choose 
$\theta \in [0, 2 \pi]$ such that $\kappa_a \prec 0$ on $\arg z = \theta$. 
A small change to $\theta$  delivers 
either $\kappa_b \prec \kappa_a$, which leads to (B), 
or $\kappa_a \prec \kappa_b \prec 0$ or $\kappa_a \prec 0 \prec \kappa_b $,  leading to (B) or (C).

Assume now that none of the $\kappa_j$ is the zero polynomial. 
Let $m^*$ be the largest of the degrees of the $\kappa_j$ and,  with no loss of generality, write 
$$
\kappa_j (z) = \alpha_j z^{m^*} + \ldots , \quad \alpha_1 \neq 0.
$$
If $\alpha_2 = 0$ then it is easy to choose a ray $\arg z = \theta $ on which 
$\kappa_1 \prec \kappa_2 \prec 0$ and, by varying $\theta$ slightly if necessary, either 
$\kappa_3 \prec 0$ or $0 \prec \kappa_3 $. 
Lemma \ref{lemkappanotsame} then implies that (A) must hold. 

Next, suppose that $\alpha_j \neq 0$ for all $j$.  
If $\alpha_2 / \alpha_1$ is not a negative real number choose a ray on which 
$\alpha_1 z^{m^*}$ and $\alpha_2 z^{m^*}$ both have negative real part and $\kappa_3 \not \simeq 0$. Shifting $\theta$ slightly gives
either $\kappa_1 \prec \kappa_2 \prec 0$ or $\kappa_2 \prec \kappa_1 \prec 0$, and Lemma \ref{lemkappanotsame} forces (A) to hold, subject to re-labelling if necessary.

Thus the proof is complete, after re-labelling if necessary, unless both $\alpha_2 / \alpha_1$ and $\alpha_3 / \alpha_1$ 
are negative real numbers, in which case the argument of the previous paragraph applies with 
$\kappa_2$ and $\kappa_3$ in place of 
$\kappa_1$ and $\kappa_2$.
\hfill$\Box$
\vspace{.1in}

\section{A decomposition of the operators $N$ and $V$}

By Lemma \ref{lemmainstep}(C), the equation $N[y]=0$ in Lemma \ref{lemgeqn}, which is satisfied by $g$, has order 
$3$, and so the asymptotics for its solutions may be complicated. However, the following
lemma gives a condition under which two linearly independent solutions of $N[y]=0$ must together solve
a second order equation, for which the asymptotics are then considerably simpler.

\begin{lem}
 \label{newdelem}
With $N$ and $V$ as in Lemma \ref{lemgeqn}, suppose that $g_1$ and $g_2$ are linearly independent (both formal or both locally analytic) solutions of $N[y]=0$ 
such that 
\begin{equation}
 \label{aa3}
g_1 V[g_2 ] - g_2 V[g_1] = d (g_2 g_1' - g_1 g_2' ), 
\end{equation}
where $d$ is rational at infinity. 
Then $g_1$ and $g_2$ solve an equation (\ref{7}) and 
\begin{equation}
 \label{a4} 
 W' + E_1 W  = 0, \quad E_1 = \frac{ \beta+ d}{\alpha} , \quad W = W(g_1, g_2 ) ,
\end{equation}
\begin{equation}
 \label{a5}
N = (D + \delta ) \circ (D^2 + E_1 D + E_0)   , 
\end{equation}
where  $E_1$,  $E_0$ and  $\delta$ are rational at infinity. If, in addition, $d$ is constant then  
\begin{equation}
 \label{a9}
V =  \alpha (D^2 + E_1 D + E_0) - d D - x E_1 - \frac{X'}{X}  ,
\quad  x = d - \frac{k-1}2  , \quad \frac{X'}X = - \frac{a_{k-1}}k . 
\end{equation}
\end{lem}
\textit{Proof.} Differentiating $W = W(g_1, g_2) = g_1 g_2' - g_1' g_2 $
gives $W' = g_1 g_2'' - g_1'' g_2 $. Thus
equation (\ref{aa3}) can be rewritten, using (\ref{5}), in the form (\ref{a4}),
with $E_1$  rational at infinity. 
Applying Lemma \ref{lem3} shows that $g_1$ and $g_2$ solve an equation (\ref{7}), 
with $E_0$ also rational at infinity. 
Because $g_1$ and $g_2$ are
linearly independent solutions of $N[y]=0$ and
(\ref{7}), the operator $N$ factorises using (\ref{Nform})
and the division algorithm for linear differential operators \cite[p.126]{Ince}
as
\begin{eqnarray}
\label{a5a}
N &=& D^3 + B_2 D^2 + B_1 D + B_0 \nonumber \\
&=& (D + \delta ) \circ (D^2 + E_1 D + E_0) \nonumber \\
&=& D^3 + (E_1+\delta) D^2 + (E_0+E_1' + \delta E_1) D + E_0' + \delta E_0   , 
\end{eqnarray}
where $\delta $ is again rational at infinity.

Now suppose that $d$ is constant. (\ref{5a}) and (\ref{5}) yield 
\begin{equation*}
 \label{a6}
h' = \alpha g''' + (\beta+\alpha') g'' + (\beta' + \gamma) g' + \gamma' g = -
\left( \frac{k-1}2 \right) \, g'' + \frac{a_{k-1}}k g' + \frac{a_{k-1}' -D_{k-2}}k g, 
\end{equation*}
where $D_{k-2}$ is rational at infinity, so that
$$
0 = \alpha g''' + \left(\beta+\alpha' + \frac{k-1}2 \right) g'' + \left(\beta' + \gamma -  \frac{a_{k-1}}k \right) g' + 
\left(\gamma' + \frac{D_{k-2} - a_{k-1}' }k \right) g .
$$
Comparing coefficients with (\ref{a5a}) leads, using Lemma \ref{lemmainstep}(C),  to
\begin{equation}
 \label{a7}
E_1+\delta = \frac{ \beta+\alpha'+(k-1)/2 }{\alpha}, \quad E_0+E_1' + \delta E_1 = \frac{\beta' + \gamma - a_{k-1}/k }{\alpha}.
\end{equation}
Now (\ref{a4}) and (\ref{a7}) deliver, with $x = d - (k-1)/2$, 
\begin{eqnarray*}
 \beta &=& \alpha E_1 - d , 
\quad 
\delta = \frac{ \beta+\alpha'+(k-1)/2 }{\alpha} - E_1 = \frac{\alpha'}{\alpha} - \frac{x}{\alpha} ,\nonumber \\
\gamma &=& \alpha ( E_0+E_1' + \delta E_1 ) - \beta' + \frac{a_{k-1}}k =
\alpha ( E_0+E_1' + \delta E_1 ) - \alpha' E_1 - \alpha E_1' + \frac{a_{k-1}}k \nonumber \\
&=& \alpha E_0 + E_1 ( \alpha' - x) - \alpha' E_1 + \frac{a_{k-1}}k = \alpha E_0 - x E_1 + \frac{a_{k-1}}k , 
\end{eqnarray*}
and the representations given here for $\beta$ and $\gamma$ yield (\ref{a9}).
\hfill$\Box$
\vspace{.1in}

\begin{lem}
 \label{d=0lemA}
Suppose that $g_1$ and $g_2$ are linearly independent (both formal or both locally  analytic) solutions of $N[y]=0$, 
such that (\ref{aa3}) holds, with $d \in \{ 0, \ldots, k-1 \}$ a constant. 
Then $g_1$ and $g_2$ solve an equation (\ref{7}), and formulas (\ref{a4}) to (\ref{a9}) hold, 
with  $E_1$,  $E_0$ and  $\delta$ rational at infinity, and $d$ satisfies $d \neq (k-1)/2$.

Suppose further that  $E_1(\infty) \neq 0$ in  (\ref{7}).
Then $g$ is given by (\ref{newestgrep}). 
If, in addition, $d=0$ or $d=k-1$, then $f$ satisfies the conclusion of Proposition \ref{mainprop}. 
\end{lem}
\textit{Proof.} Lemma \ref{newdelem} gives
the equation (\ref{7}) solved by $g_1$ and $g_2$, as well as  formulas  (\ref{a4}) to (\ref{a9}). 
Since $N[g]=0$, but $g'' + E_1g' + E_0 g \not \equiv 0$, (\ref{5}),  (\ref{a5}) and (\ref{a9}) deliver
$$
g'' + E_1g' + E_0 g = e^\tau  , \quad \tau' = - \delta    , 
$$
and 
$$ 
- \frac{f'}{f} = \frac{h}{g} = \frac{V[g]}g = \frac{\alpha e^\tau}g - \frac{dg'}{g} -x E_1 + \frac{a_{k-1}}k =
\frac{\alpha e^\tau}g - \frac{dg'}{g} + a^*  .
$$
Here the function $a^* = -x E_1 + a_{k-1}/k$  is rational at infinity and $-f'/f +dg'/g-a^*$ continues without zeros in some 
$\Omega(r_2)$. Hence Lemma \ref{propf'not0} and (\ref{a9}) imply that $d \neq (k-1)/2$ and $x \neq 0$.
Finally, 
if $E_1(\infty) \neq 0$ in  (\ref{7})
then  $a^*(\infty) \neq 0$ and the remaining assertions of Lemma \ref{d=0lemA} follow from Lemma \ref{propf'not0}.
\hfill$\Box$
\vspace{.1in}

\section{Analytic solutions decaying in the same sector }\label{getsecondorder}

This section determines conditions under which Lemma \ref{d=0lemA} may be applied
with analytic solutions of $N[y]=0$.

\begin{lem}
 \label{lemtwogsmall}
Assume 
that there exist linearly independent analytic solutions $g_1$, $g_2$ of $N[y] = 0$, such that 
both tend to $0$ transcendentally fast as $z \to \infty$ in the same sector $S$. 
Then $g_1$ and $g_2$ satisfy the hypotheses of Lemma \ref{d=0lemA}, for some constant $d \in \{ 0, \ldots, k-1 \}$, and  
solve an equation (\ref{7}),  
with  $E_1$ and  $E_0$  rational at infinity and $E_1(\infty) \neq 0$. 
Moreover, formulas (\ref{a4}) to (\ref{a9}) hold, and $d$ and $g$ satisfy $d \neq (k-1)/2$ and  (\ref{newestgrep}). Finally,  
if $d=0$ or $d=k-1$ then  the conclusion of Proposition~\ref{mainprop} holds. 
\end{lem}
In the context of Section \ref{expparts}, Lemma \ref{lemtwogsmall} applies if 
there is a repeated non-trivial exponential part among the $\kappa_j$, 
or if (\ref{different1}) holds for some ray $\arg z = \theta$. 
\\\\
\textit{Proof.} 
Choose $z_0 \in S$ such that $z_0$ is not 
a singular point for any of the operators $L, M, N$, and such that
\begin{equation}
 \label{aa1}
g_1(z_0) g_2(z_0) \neq 0, \quad W(g_1, g_2) (z_0) \neq 0. 
\end{equation}
Let $w $
lie close to $z_0$. Then $g_w(z) = g_2(w) g_1(z) - g_1(w) g_2(z)$ tends to $0$ transcendentally fast as $z \to \infty$ in $S$
and, by Lemma \ref{lem1}, $g_w$ annihilates a  solution
$f_w \not \equiv 0$ of $L[y]=0$, with
\begin{equation}
 \label{aa2}
\frac{f_w'(z)}{f_w(z)} = - \frac{V[g_w](z)}{g_w(z)} = \frac{ g_1(w) V[g_2 ](z) - g_2(w) V[g_1] (z)}{ g_2(w) g_1 (z) - g_1(w) g_2(z)} .  
\end{equation}
Let
$$
G_0(z)= \frac{ g_1 (z)V[g_2 ](z) - g_2 (z)V[g_1](z) }{ g_2(z) g_1'(z) - g_1 (z)g_2'(z)} .
$$
The second condition of (\ref{aa1}) implies that $w$ is a simple zero of $g_w$, and by (\ref{aa2}) the 
residue of $f_w'/f_w$ at $w$
is $G_0(w)$, which must belong to the set $\{ 0, \ldots, k-1 \}$. Since this holds for all $w$ near $z_0$,
the function $G_0$ is a constant $d \in \{ 0, \ldots, k-1 \}$, and so $g_1$ and $g_2$ satisfy the hypotheses of Lemma~\ref{d=0lemA},
and hence solve an equation (\ref{7}). Cauchy's estimate for derivatives shows that
$W(g_1, g_2)$ tends to $0$ transcendentally fast in a subsector of $S$, 
which gives $E_1(\infty) \neq 0$ by Abel's identity.
The remaining assertions hold by Lemma \ref{d=0lemA}.
\hfill$\Box$
\vspace{.1in}

\section{Decaying solutions with different exponential parts}
\label{twosmall}

This section will deal with one case of the situation in Section \ref{getsecondorder}, 
in which two linearly independent solutions of $N[y]=0$ decay in the same sector and
have different exponential parts, corresponding to (\ref{different1}) in Lemma \ref{lemkappaalldifferent}.
The case of a repeated non-trivial exponential part will be addressed in Section \ref{repeatnontriv}.
The methods of this section are heavily influenced by \cite{Bru}, but a decisive role will be played by Lemma \ref{lemtwogsmall}
and the second order equation (\ref{7}).

\begin{prop}
 \label{prop3}
Assume
that there exists
a ray $\arg z = \theta $ on which the exponential parts 
$\kappa_j$  for the equation $N[y] =0$ satisfy (\ref{different1}). 
Then 
$f$ satisfies the conclusions of Proposition~\ref{mainprop}. 
\end{prop}

To prove Proposition \ref{prop3}, note first that if $\theta$ is varied slightly, then (\ref{different1})
continues to hold.
Take canonical formal solutions $g_1, g_2$ of $N[y]=0$ with exponential parts $\kappa_1, \kappa_2$. 
By (\ref{different1}) the exponential parts for $N[y]=0$ are pairwise distinct, and there exist linearly independent 
analytic solutions $G_1, G_2$ of $N[y]=0$ which are asymptotic to $g_1$ and $g_2$ respectively 
on a sector centred on the ray $\arg z = \theta $, and so tend to $0$ transcendentally
fast there, by (\ref{different1}). Thus the hypotheses of 
Lemma \ref{lemtwogsmall} are satisfied, and therefore so are those of Lemma \ref{d=0lemA},
for some $d $ in $\{ 0, \ldots, k-1 \}$, which gives rise to an equation (\ref{7}) satisfied by the $G_j$. Computing series
representations for $0 = G_j''+E_1 G_j'+E_0 G_j $ shows that the $g_j$ also solve (\ref{7}). 
Thus the exponential parts 
$\kappa_1$ and $\kappa_2$ correspond to the equation (\ref{7}), 
while $\kappa_3$ is from the third canonical formal solution of $N[y]=0$.
Furthermore, $V$ satisfies (\ref{a9}).

Next, let the operators $L$, $M$ have canonical formal solutions
with exponential parts $q_j$, $s_j$ respectively,  labelled so that 
\begin{equation}
 q_1 \preceq q_2 \preceq \ldots \preceq q_k , 
\quad s_1 \preceq s_2 \preceq \ldots \preceq s_k  ,
\label{yz2}
\end{equation}
on $\arg z = \theta $ (the last phrase will be omitted henceforth).
The $q_j$ and $s_j$ are polynomials in $z$ with zero constant term. 
It may be assumed that $\theta$ is chosen so that if 
$\widetilde p_1, \widetilde p_2 \in \{ q_1, \ldots, q_k, s_1, \ldots, s_k \}$ 
and $ \widetilde p_1 - \widetilde p_2 \not \equiv 0 $ then 
$\widetilde p_1 \prec \widetilde p_2$  or $\widetilde p_2 \prec \widetilde p_1$.

\begin{lem}
 \label{lemflambdafmu}
There exists $\lambda \in \{ 1, \ldots, k \}$ such that the canonical formal solution 
$g_1$ of (\ref{7}) with exponential part $\kappa_1$ annihilates a  canonical
formal solution $f_\lambda$ of $L[y] = 0$ with
exponential part $q_\lambda$, 
and the exponential parts 
for $M[y]=0$ are
\begin{equation}
 \label{yz3}
q_j + \kappa_1 \quad (j \neq \lambda), \quad q_\lambda - (k-1) \kappa_1 .
\end{equation}
Moreover, this $f_\lambda$ may be assumed to be $g_1^d W^{-x}X$, where $W$, $x$ and $X$ are as in Lemma \ref{newdelem}. 
Furthermore, there exists $\mu \in \{ 1, \ldots, k \}$ such that the canonical formal solution  
$g_2$ of (\ref{7}) with exponential part $\kappa_2$ annihilates a canonical
formal solution $f_\mu = g_2^d W^{-x}X$ of $L[y] = 0$, 
with exponential part $q_\mu$, 
while the exponential parts for $M[y]=0$ are
\begin{equation}
 \label{yz4}
q_j + \kappa_2 \quad (j \neq \mu), \quad q_\mu - (k-1) \kappa_2 .
\end{equation}
\end{lem}
\textit{Proof.} Since the exponential parts for $N[y]=0$ are pairwise distinct, 
$g_1$ and $g_2$ both have non-zero exponential parts and are free of logarithms. 
Thus  Lemma \ref{lemexppart} gives a  canonical
formal solution $f_\lambda$ of $L[y] = 0$, with
exponential part $q_\lambda$, such that $g_1$ annihilates $f_\lambda$ and the exponential parts 
for $M[y]=0$ are given by (\ref{yz3}). Moreover, since $g_1$ is a solution of (\ref{7}), solving $0 = f_\lambda' g_1 + f_\lambda V[g_1]$ in the light of 
(\ref{a9}) shows that $f_\lambda$ is a constant multiple of $g_1^d W^{-x}X$.
The same argument works for $g_2$ and $f_\mu$.
\hfill$\Box$
\vspace{.1in}

\begin{lem}
\label{lemdr1} 
The integer  $\lambda $ is $ 1$. 
\end{lem}
\textit{Proof.} 
Suppose not:
then an exponential part $q_1 + \kappa_1$ occurs in the list (\ref{yz3}). But this term,  in view of (\ref{different1}) and
(\ref{yz2}),
cannot be realised as $q_j + \kappa_2 $ or $ q_\mu - (k-1) \kappa_2 $. 
\hfill$\Box$
\vspace{.1in}

\begin{lem}
\label{lemdr2} 
The $q_j$ satisfy $q_j + \kappa_1 \preceq q_1 - (k-1) \kappa_1 $ for $2 \leq j \leq k$. 
\end{lem}
\textit{Proof.}
Suppose that this is not the case. Then the term $q_k + \kappa_1$, which does occur in the  list (\ref{yz3}), must be maximal
according to the ordering $\preceq$.
But  (\ref{different1}) implies that 
\begin{equation}
 \label{qk+kappa1}
q_k + \kappa_1 \prec  q_k + \kappa_2 \prec q_k - (k-1) \kappa_2 .
\end{equation}
This is a contradiction since  the second or third term in (\ref{qk+kappa1}) occurs in the list (\ref{yz4}). 
\hfill$\Box$
\vspace{.1in}

Thus by (\ref{yz3}) the $s_j$ in  (\ref{yz2}) can now be written as 
\begin{equation}
 \label{yz5}
s_1 = q_2 + \kappa_1,  \quad \ldots, \quad s_{k-1} = q_k + \kappa_1, \quad s_k = q_1 - (k-1) \kappa_1. 
\end{equation}
Note that each of these relations initially holds with $\simeq $ in place of $=$, but may be assumed to be an identity,
by the remark following (\ref{yz2}). The same property will subsequently
be used on a number of occasions without explicit reference. 

\begin{lem}
\label{lemdr3} 
The exponential part $s_\mu$ satisfies $s_\mu = q_\mu - (k-1) \kappa_2 $. 
\end{lem}
\textit{Proof.} 
Suppose first that $ q_\mu - (k-1) \kappa_2 \prec s_\mu$. Then 
(\ref{different1}) and (\ref{yz2})  give $\mu > 1$ and 
$$
q_{1} + \kappa_2 \preceq \ldots \preceq q_{\mu-1}  + \kappa_2 \prec q_\mu - (k-1) \kappa_2 \prec s_\mu ,
$$
in which all of the first $\mu$ terms occur in the list (\ref{yz4}). 
Hence the second list in (\ref{yz2}) 
includes $ \mu$ terms $\widetilde s$ all  satisfying $ \widetilde s \prec s_\mu $, which is a contradiction. 

Now suppose that $s_\mu \prec q_\mu - (k-1) \kappa_2  $.
Then $\mu < k$ and in the list 
(\ref{yz4}) there are at least $ \mu $ terms $\widetilde s$
all satisfying $ \widetilde s \preceq s_\mu \prec q_\mu - (k-1) \kappa_2 $. 
Of these, $\mu-1$ are $q_{1} + \kappa_2, \ldots, q_{\mu-1} + \kappa_2$ (this list being void if $\mu=1$), and 
it must be the case that $q_{\mu+1} + \kappa_2 \preceq s_\mu
\prec q_\mu - (k-1) \kappa_2  $. But then (\ref{different1}) and (\ref{yz5}) yield a contradiction via
$$
s_\mu = q_{\mu+1} + \kappa_1 \prec q_{\mu+1} + \kappa_2 \preceq s_\mu .
$$
\hfill$\Box$
\vspace{.1in}

Lemma \ref{lemdr3} implies that among the  $q_j + \kappa_2$ ($j \neq \mu$) there are at least 
$\mu-1$ terms $\widetilde s$ 
with 
$\widetilde s \preceq q_\mu - (k-1) \kappa_2 $, and if $\mu>1$ these must  include $q_1 + \kappa_2, \ldots, q_{\mu-1} + \kappa_2$; similarly, there are at least $k - \mu$ terms
with
$q_\mu - (k-1) \kappa_2  \preceq \widetilde s$, and if $\mu < k$ these must include $q_{\mu+1} + \kappa_2, \ldots, q_k + \kappa_2$. It follows that  
\begin{equation}
 \label{yz8}
s_j = q_j + \kappa_2 \quad (j \neq \mu), \quad s_\mu = q_\mu - (k-1) \kappa_2 .
\end{equation}

\begin{lem}
 \label{lemdr44}
The integers $d$ and $\mu$ are related by  $d= \mu-1$. 
\end{lem}
\textit{Proof.} By Lemmas \ref{lemflambdafmu} and \ref{lemdr1} the canonical formal solutions $f_\lambda$ and $f_\mu$  of $L[y]=0$
annihilated by $g_1$ and $g_2$  have exponential parts $q_1$ and $q_\mu$ respectively. The quotient 
$f_\mu/f_\lambda = ( g_2/g_1 )^{d}$ has exponential part $d( \kappa_2 - \kappa_1 )$, which implies that $d( \kappa_2 - \kappa_1 ) = q_\mu - q_1$.
If $\mu=1$ this gives $d=0$ since  $\kappa_1 \neq \kappa_2$.
For $\mu >1$,  (\ref{yz5}) and (\ref{yz8}) yield
$$
s_1 = q_2 + \kappa_1,  \quad \ldots, \quad s_{\mu-1} = q_\mu  + \kappa_1, \quad
s_1 = q_1 + \kappa_2 , \quad \ldots, \quad s_{\mu-1} = q_{\mu-1} + \kappa_2,
$$
and so 
$$
\kappa_2 - \kappa_1 = q_2 - q_1 = \ldots = q_\mu - q_{\mu-1} , \quad 
d( \kappa_2 - \kappa_1 ) = q_\mu - q_1 = (\mu-1) (\kappa_2 - \kappa_1 ). 
$$
\hfill$\Box$
\vspace{.1in}

By (\ref{different1}) there exists a canonical formal solution $g_3$ of $N[y]=0$ which is free of logarithms
and has  exponential part $\kappa_3$. 

\begin{lem}  \label{lemdr5a}
The exponential part $\kappa_3$ is not the zero polynomial, and case (A) applies in (\ref{different1}).
\end{lem}
\textit{Proof.} 
Suppose that $\kappa_3 \equiv 0$. Then Lemma \ref{lemexppart} and (\ref{yz8}) give at least one $j$ with 
$q_j = s_j = q_j + \kappa_2 \prec q_j$, a contradiction. 
\hfill$\Box$
\vspace{.1in}

By Lemmas \ref{lemexppart} 
and \ref{lemdr5a}, there exists $\nu$ such that $g_3$ 
annihilates a canonical formal solution of $L[y]=0$ with exponential part $q_\nu$, and 
the exponential parts for $M[y]=0$ are
\begin{equation}
 \label{yz9}
q_j + \kappa_3 \quad (j \neq \nu), \quad q_\nu - (k-1) \kappa_3 .
\end{equation}

\begin{lem}
 \label{lemnewbrug}
Assume that 
$2 \leq \mu \leq k-1$. Then  $\nu = k$  and $s_1 = q_k-(k-1) \kappa_3$.
\end{lem}
\textit{Proof.} Suppose first that $\nu < k$. Then the list (\ref{yz9}) includes $q_k + \kappa_3$, which must be maximal 
with respect to the ordering $\preceq$, since $0 \prec \kappa_3$. But $\mu \neq  k$ by assumption, which  gives
$$
q_k + \kappa_3 = s_k = q_k + \kappa_2
$$
using (\ref{yz8}), and this contradicts (\ref{different1}). Thus $\nu = k$ in (\ref{yz9}). 

Now suppose that $s_1 \neq q_k-(k-1) \kappa_3$. Then  $s_1 \prec q_k-(k-1) \kappa_3$ and 
so $s_1 = q_1 + \kappa_3$, whereas  (\ref{yz8}) gives $s_1 = q_1 + \kappa_2$ since $\mu \neq  1$, again contradicting (\ref{different1}).
\hfill$\Box$
\vspace{.1in}

\begin{lem}
 \label{lemmuk-1}
If $k \geq 4$ then $\mu=1$ or $\mu=k$.  
\end{lem}
\textit{Proof.}
Suppose instead that $2 \leq \mu \leq k-1$. Then, by Lemma \ref{lemnewbrug}, the list (\ref{yz9}) consists of
\begin{equation}
 \label{yz12a}
s_1 = q_k - (k-1) \kappa_3 , \quad s_2 = q_1 + \kappa_3, \quad \ldots, \quad s_{k} = q_{k-1} + \kappa_3 . 
\end{equation}
Using (\ref{yz5}), (\ref{yz8}) and (\ref{yz12a}) gives
\begin{equation}
 \label{yz15a}
s_{\mu+1} = q_\mu + \kappa_3 = q_{\mu+1} + \kappa_2  , \quad q_{\mu} - q_{\mu+1} =  \kappa_2 -  \kappa_3, 
\end{equation}
and 
\begin{equation}
 \label{yz16a}
s_\mu = q_{\mu+1} + \kappa_1 = q_\mu - (k-1) \kappa_2, \quad \quad q_{\mu} - q_{\mu+1} = \kappa_1 + (k-1) \kappa_2.
\end{equation}
Define $\tau$ as follows: if $2 \leq \mu \leq k-2$  take $\tau = k-1$, and if $\mu = k-1$  choose $\tau=1$. 
In either case $\mu \neq \tau, \tau+1 $, since $k \geq 4$ by assumption.
Thus (\ref{yz5}), 
(\ref{yz8}) and (\ref{yz12a}) deliver 
\begin{equation}
 \label{yz19}
s_{\tau+1} = q_{\tau+1} + \kappa_2 = q_\tau+\kappa_3 , \quad 
q_{\tau+1} - q_\tau = \kappa_3 - \kappa_2 , 
\end{equation}
in addition to
\begin{equation}
 \label{yz20}
s_{\tau} = q_{\tau} + \kappa_2 = q_{\tau+1} + \kappa_1 , \quad 
q_{\tau+1} - q_\tau = \kappa_2 - \kappa_1 . 
\end{equation}
Combining (\ref{yz15a}),  (\ref{yz16a}), (\ref{yz19}) and (\ref{yz20}) yields 
$$
\kappa_2 -  \kappa_1   =  \kappa_3 -  \kappa_2  = - \kappa_1 - (k-1) \kappa_2, 
$$
contradicting the fact that $ \kappa_2 \prec 0$.
\hfill$\Box$
\vspace{.1in}

Thus $d$ must be $0$ or $k-1$: this  follows from  
Lemmas \ref{lemdr44} and \ref{lemmuk-1} when $k \geq 4$, while if $k=3$ then  
Lemma \ref{lemtwogsmall} forces $d \neq (k-1)/2 = 1$.  
Hence
the conclusion of Proposition \ref{mainprop} holds by Lemma \ref{lemtwogsmall} and  the proof of Proposition \ref{prop3} is complete.

\hfill$\Box$
\vspace{.1in}

\section{The case of a repeated non-trivial exponential part}\label{repeatnontriv}

Suppose that $\kappa$ is a repeated non-trivial exponential part for the equation $N[y]=0$. 
Then it is possible to choose
a ray $\arg z = \theta \in \R$ on which $\kappa \prec 0$, and linearly independent analytic solutions
$g_1, g_2$ of $N[y]=0$, each with exponential part $\kappa$ near $\arg z = \theta $. 
It then follows from Lemma \ref{lemtwogsmall}
that 
$g$ is given by (\ref{newestgrep}), which yields  
$$
0 = N[g] = \beta_1 H_1 e^{\omega_1 P} + \beta_2 H_2 e^{\omega_2 P} + \beta_3 H_3 e^{\omega_3 P}  , \quad \beta_j , \omega_j \in \C^*, 
$$
in which   the $\omega_j$ are pairwise distinct, while $P'$ and the $H_j$ are rational at infinity and 
$ H_je^{\omega_j P} = N[e^{\omega_j P} /P' ]$.  This forces each $H_j$ to vanish identically, so that the equation $N[y]=0$ has three pairwise distinct exponential
parts for its solutions, which is a contradiction. 
\hfill$\Box$
\vspace{.1in}

\section{Two lemmas concerning trivial exponential parts}\label{onetrivial}

If at least one of the three exponential parts arising from the equation $N[y]=0$ is trivial
(that is, the zero polynomial), then it is not necessarily the case that $N[y]=0$ will
have two linearly independent solutions decaying in the same sector, so that a second order
equation (\ref{7}) may not be available. The approach to this case will combine  Lemma \ref{lemintegervalued}
with some ideas  from \cite{Bru}.

\begin{lem}
 \label{prop2a}
Assume that two exponential parts 
$\kappa_1, \kappa_2$ arising from the equation $N[y] =0$ are such that
$\kappa_2 $ is the zero polynomial, while 
\begin{equation}
 \label{zz1a}
\kappa_1 \prec 0  \quad \hbox{or} \quad 0
\prec \kappa_1
\end{equation}
on a ray $\arg z = \theta $.
Let the operators $L$, $M$ have canonical formal solutions
with exponential parts  as in (\ref{yz2}). Then 
the exponential parts  for $M$ are as in (\ref{yz3}), while  
\begin{equation}
 \label{triv1}
\hbox{$s_j = q_j$ for each $j$ }
\end{equation}
and the following additional conclusions hold.

If $\kappa_1 \prec 0$ in (\ref{zz1a}) then  $\lambda = 1$ and 
\begin{equation}
 \label{alldifferent}
q_1 = s_1 = q_2 + \kappa_1, \quad 
\ldots,
\quad 
q_{k-1} = s_{k-1} = q_k + \kappa_1, \quad q_k = s_k =   q_1 - (k-1) \kappa_1 .
\end{equation}
If $0 \prec \kappa_1 $ in (\ref{zz1a}) then  $\lambda = k$ and 
\begin{equation}
 \label{alldifferent2}
q_1 = s_1 = q_k - (k-1) \kappa_1 , \quad 
q_2 = s_2 = q_1 + \kappa_1, 
\ldots,
\quad 
q_{k} = s_{k} = q_{k-1} + \kappa_1 .
\end{equation}
\end{lem}
\textit{Proof.} 
First observe that $N[y]=0$ has two canonical formal solutions which are free of logarithms and have exponential parts $0$ and $\kappa_1$ respectively.
Thus (\ref{yz3}) and (\ref{triv1}) hold by Lemma \ref{lemexppart}.
Assume  that $\kappa_1 \prec 0$ in (\ref{zz1a}). If $\lambda \neq 1$ then an 
exponential part $q_1 + \kappa_1$ occurs in the list (\ref{yz3}),
but this term,  in view of (\ref{zz1a}),
cannot be realised as $q_j  $ for any $j$, contradicting (\ref{triv1}). 
Now suppose that $s_k \neq q_1 - (k-1) \kappa_1 $; then
$s_k = q_k + \kappa_1  $, again contradicting (\ref{triv1}). 

Now assume  that $0 \prec \kappa_1 $ in (\ref{zz1a}). Then
$\lambda$ must be $k$, since otherwise an exponential part $q_k + \kappa_1$ occurs in (\ref{yz3}), 
contradicting (\ref{triv1}).  
Moreover, $q_1 = s_1 = q_k - (k-1) \kappa_1$, because the contrary case forces
$s_1 = q_1 + \kappa_1$, which again contradicts (\ref{triv1}). 
\hfill$\Box$
\vspace{.1in}

\begin{lem}
 \label{prop2}
If there exists  a ray $\arg z = \theta $ on which the  three exponential parts 
arising from the equation $N[y] =0$ satisfy (\ref{different2}), 
then  $\kappa_3 = - \kappa_1$.  
\end{lem}
\textit{Proof.}  Assuming the existence of such a ray, 
let the operators $L$, $M$ have 
exponential parts  as in (\ref{yz2}). 
Now (\ref{alldifferent}) and (\ref{alldifferent2}) yield
$$
q_1 = s_1 = q_k - (k-1) \kappa_3, \quad q_k = s_k = q_1 - (k-1) \kappa_1, \quad \kappa_3 = - \kappa_1 .
$$
\hfill$\Box$
\vspace{.1in}

\section{The case where (\ref{different2}) holds}\label{onetrivial3}

This section will deal with the case where there exists a ray for which conclusion (\ref{different2}) arises in Lemma
\ref{lemkappaalldifferent}. In this situation Lemma \ref{prop2} makes it possible 
to assume that the exponential parts for $N[y]=0$ are  $P$, $0$ and $-P$, where $P$ is a polynomial in $z$ of 
positive degree $\rho$. Hence  $N[y] = 0$ has canonical formal solutions which are free of logarithms and satisfy
\begin{equation}
\label{ab2}
 u_1 (z) = z^{\eta_1} e^{P(z)} (1+ \ldots ), \quad 
u_2 (z) = z^{\eta_2}  (1+ \ldots ), \quad 
u_3 (z) = z^{\eta_3} e^{-P(z)} (1+ \ldots ). 
\end{equation}
Since $N[g]=0$, the order of growth of $g^k = f/F$ is $\rho\left(g^k\right) = \rho$. 
Choose a ray $\arg z = \theta_0$ on which ${\rm Re} \, P(z) = O( |z|^{\rho-1})$ as $|z| \to \infty$,
such that $f$ has a sequence of  poles (and so $g$ has a sequence of simple zeros) tending to $\infty$ in the sector $| \arg z - \theta_0 | \leq \pi / 2 \rho $.
Take a sector $\Sigma$ given by $| \arg z - \theta_0 | \leq \pi /  \rho - \delta_1 $, where $\delta_1$ is small and positive,
and write 
\begin{equation}
 \label{ac1}
g = U_1 + U_2 + U_3 , \quad U_j = b_j \phi_j, \quad b_j \in \C, \quad \phi_j \sim   u_j ,
\end{equation}
in which 
the $\phi_j$ are analytic solutions on $\Sigma$, and the last relation holds in the sense of asymptotic series,
as in Section \ref{wasowthm}. 
Here the fact that the asymptotics for $N[y]=0$ may be extended to hold in $\Sigma$ follows
from the work of Jurkat \cite{jurkat}: in the present case, where the exponential parts are $P$, $0$
and $-P$, it is relatively simple to establish, using the Phragm\'en-Lindel\"of principle.
Since $g$ has infinitely many zeros in $\Sigma$, at least two of the $b_j$,
and so at least one of $b_1$ and $b_3$, must be non-zero. By replacing $P$ by $-P$, it may be assumed that $b_1 \neq 0$. 

Now take a ray $\arg z = \theta$ lying in $\Sigma$, on which $P \prec  0 \prec  -P$, and apply Lemma \ref{prop2a} with $\kappa_1 = P$.
It follows from (\ref{triv1}) and (\ref{alldifferent})  that 
\begin{equation}
 \label{ab5}
q_2 = q_1 -P, \quad q_3 = q_2-P = q_1-2P, \quad \ldots, \quad q_k = q_{k-1} -P = q_1 - (k-1)P , \quad s_j = q_j , 
\end{equation}
and so, by (\ref{abel}),  
\begin{equation}
 \label{ab6}
0 = q_1 + \ldots + q_k = k q_1 - \left( \frac{k(k-1)}2 \right) P, \quad q_1 = \left( \frac{k-1}2 \right) P , \quad q_k = - \left( \frac{k-1}2 \right) P.
\end{equation}
Hence  the equations $L[y]=0$, $M[y]=0$ have canonical formal solutions 
\begin{equation}
 \label{ab3}
f_j(z) = z^{\lambda_j} e^{q_j(z)} (1 + \ldots ), \quad 
w_j(z) = z^{\mu_j} e^{q_j(z)} (1 + \ldots ),
\end{equation}
respectively, 
in which the $q_j$ are pairwise distinct. Since $a_{k-1} = A_{k-1}$,
Lemma \ref{wronskianlem3} implies that 
\begin{equation}
 \label{ab4}
\lambda_1 + \ldots + \lambda_k = \mu_1 + \ldots +\mu_k .
\end{equation}
Write $v_j = V[u_j]$. By Lemmas \ref{lemexppart} and \ref{prop2a}, $u_1$, $u_3$ annihilate $f_1$, $f_k$ respectively,
and (\ref{ab6}) gives 
\begin{equation}
 \label{ab7}
\frac{v_1(z)}{u_1(z)} = - \frac{f_1'(z)}{f_1(z)} = \widehat c_1 z^{\rho-1} + \ldots , \quad 
\frac{v_3(z)}{u_3(z)} = - \frac{f_k'(z)}{f_k(z)} = \widehat c_3 z^{\rho-1} + \ldots , 
\end{equation}
where $ \widehat c_1$, $ \widehat c_3$ are non-zero constants. 
It follows from 
(\ref{5}) and (\ref{ab2}) that $v_2/u_2$ is given by a (possibly vanishing) formal series in descending integer powers of $z$
of the form 
\begin{equation}
 \label{ab10}
\frac{v_2(z)}{u_2(z)} = \alpha (z) \frac{u_2''(z)}{u_2(z)} + \beta(z) \frac{u_2'(z)}{u_2(z)} + \gamma (z) = 
c_N z^N + \ldots . 
\end{equation}

\begin{lem}
 \label{lemk=3}
The integer $k$ is at least $4$.
\end{lem}
\textit{Proof.} Suppose that $k=3$; then  (\ref{ab5}) and (\ref{ab6}) lead to  
\begin{equation*}
 \label{ab15}
q_1 = P, \quad q_2 = 0, \quad q_3 = -P .
\end{equation*}
Now write, using (\ref{ab2}), (\ref{ab3}) and (\ref{ab7}),
$$
f_3'(z) u_1(z) + f_3(z) v_1(z) = f_3(z) u_1(z) \left( \frac{f_3'(z)}{f_3(z)} - \frac{f_1'(z)}{f_1(z)} \right) 
=  z^{\lambda_3 + \eta_1 } (1 + \ldots ) \left( \frac{f_3'(z)}{f_3(z)} - \frac{f_1'(z)}{f_1(z)} \right) 
$$
and 
$$
f_1'(z) u_3(z) + f_1(z) v_3(z) = f_1(z) u_3(z) \left( \frac{f_1'(z)}{f_1(z)} - \frac{f_3'(z)}{f_3(z)} \right) 
=  z^{\lambda_1 + \eta_3 } (1 + \ldots ) \left( \frac{f_1'(z)}{f_1(z)} - \frac{f_3'(z)}{f_3(z)} \right) .
$$
Each of these is a formal solution of $M[y]=0$, with zero exponential part, and so a constant multiple of $w_2$.
But this implies that  
$\lambda_1 + \eta_3 = \lambda_3 + \eta_1$ and 
$f_3'u_1 + f_3 v_1 = - (f_1'u_3 + f_1v_3)$,
so that 
$f_1u_3 = f_3u_1$,
which  leads in turn to 
$$
\frac{V[u_3]}{u_3} - \frac{V[u_1]}{u_1} = 
\frac{f_1'}{f_1} - \frac{f_3'}{f_3} = \frac{u_1'}{u_1} - \frac{u_3'}{u_3},
\quad
u_1 V[u_3] - u_3 V[u_1] = u_3 u_1' - u_1u_3'.
$$
Hence (\ref{aa3}) holds, with $g_1 = u_1$, $g_2 = u_3$ and $d=1$. Thus the hypotheses of 
Lemma \ref{d=0lemA} are satisfied, so that $d \neq (k-1)/2 = 1$, a contradiction.
\hfill$\Box$
\vspace{.1in}

\begin{lem}
 \label{sofar}
One of the following two conclusions holds, in which $\rho = \deg P > 0$:\\
(A) $\eta_2 = -N \leq - \rho$ and 
$\eta_1 + \eta_3 = - 2 (\rho-1) $;\\
(B) $\eta_1 + \eta_3 - 2 \eta_2 =0 $ and $f$ has order of growth $\rho$. 
\end{lem}
\textit{Proof.} (\ref{ab5}) and (\ref{ab7}) show that, for $j=2, \ldots, k$, the term
$$
f_j' u_1 + f_j v_1 = f_j u_1 \left( \frac{f_j'}{f_j} - \frac{f_1'}{f_1} \right) 
$$
is a canonical formal solution of $M[y]=0$ with exponential part $q_j + P = q_{j-1}$, and so is a constant multiple of $w_{j-1}$.
This delivers, using (\ref{ab5}), (\ref{ab3}) and (\ref{ab4}), 
\begin{equation}
 \label{ab8}
\mu_1 = \lambda_2 + \eta_1 + \rho-1 , \quad \ldots, \quad 
\mu_{k-1} = \lambda_k + \eta_1 + \rho-1 , \quad 
\mu_k = 
\lambda_1 - (k-1) ( \eta_1 + \rho-1) . 
\end{equation}
In the same way, for $j=1, \ldots, k-1$, the term $f_j'u_3 + f_j v_3$ has exponential part $q_j -P = q_{j+1}$, and so is a constant multiple of $w_{j+1}$, which 
yields
\begin{equation}
 \label{ab9}
\mu_2 = \lambda_1 + \eta_3 + \rho-1 , \quad \ldots, \quad 
\mu_{k} = \lambda_{k-1} + \eta_3 + \rho-1 , \quad 
\mu_1 = 
\lambda_k - (k-1) ( \eta_3 + \rho-1) . 
\end{equation}

Suppose first that $N \geq \rho$ and $c_N \neq 0$ in (\ref{ab10}).
In this case (\ref{ab5}) and (\ref{ab3}) show that
$u_2$ cannot annihilate any of the  $f_j$, and that each $f_j' u_2 + f_j v_2$  is a canonical formal solution of $M[y]=0$ with exponential part~$q_j$,
and so a constant multiple of $w_j$.
This implies in view of (\ref{ab4}) that 
\begin{equation*}
 \label{ab11}
\mu_j = \lambda_j + \eta_2 + N \quad (j=1, \ldots, k), \quad \eta_2 = -N. 
\end{equation*}
Moreover, (\ref{ab8}) and (\ref{ab9}) now lead to 
$$
\lambda_k = \mu_k = \lambda_1 - (k-1) ( \eta_1 + \rho-1) , \quad 
\lambda_1 = \mu_1 = \lambda_k - (k-1) ( \eta_3 + \rho-1) ,
\quad \eta_1 + \rho-1 = - ( \eta_3 + \rho-1 ),
$$
so that $\eta_1 + \eta_3 = - 2 ( \rho -1)$ and conclusion (A) holds.

Now 
suppose that $N \leq \rho -1 $  in (\ref{ab10}): this case will lead to conclusion (B),
and encompasses the 
possibility that $v_2/u_2$ vanishes identically. The first step is to show that
the order of growth of $f$ is $\rho$. 
Since $g^k$ has order $\rho$ it follows from (\ref{2}) that the order of $f$ is at least $\rho$. 
It suffices to show that in (\ref{5}) the coefficients (which are rational at infinity) satisfy 
\begin{equation}
 \label{ab13}
\alpha(z) = O( |z|^{1-\rho} ), \quad \beta (z) = O(1), \quad 
\gamma(z) = O( |z|^{\rho-1}) \quad \hbox{as $z \to \infty$}, 
\end{equation}
because if this can be established then $\rho(f) \leq \rho$ follows from (\ref{5}), the Wiman-Valiron theory \cite{Hay5} applied to $1/f$,  
and standard estimates \cite{Gun} for logarithmic derivatives of $g^k$ and $g$.  

To prove (\ref{ab13}) 
use (\ref{ab2}), (\ref{ab7}) and (\ref{ab10}) to write 
\begin{eqnarray*}
\alpha (z) P'(z)^2 (1+O(z^{-1})) + \beta (z) P'(z) (1+O(z^{-1})) + \gamma (z) &=& \frac{v_1(z)}{u_1(z)} = O( z^{\rho-1} ) , \\
\alpha (z) P'(z)^2 (1+O(z^{-1})) - \beta (z) P'(z) (1+O(z^{-1})) + \gamma (z) &=& \frac{v_3(z)}{u_3(z)} = O( z^{\rho-1} ) , \\
\alpha (z) O( z^{-2} )  + \beta (z)O(z^{-1}) + \gamma (z) &=& \frac{v_2(z)}{u_2(z)} = O( z^N ) =  O( z^{\rho-1} ) . 
\end{eqnarray*}
Here $O( z^\omega )$ denotes any formal series in descending integer powers of $z$ with leading power at most $\omega \in \Z$.
Eliminating $\gamma$ via the last equation yields 
 \begin{eqnarray*}
\alpha (z) P'(z)^2 (1+O(z^{-1})) + \beta (z) P'(z) (1+O(z^{-1}))  &=& O( z^{\rho-1} ) , \\
\alpha (z) P'(z)^2 (1+O(z^{-1})) - \beta (z) P'(z) (1+O(z^{-1}))  &=&  O( z^{\rho-1} ) , 
\end{eqnarray*}
and now (\ref{ab13}) follows from Cramer's rule.

Next,  since $N \leq \rho-1$, 
(\ref{ab5}) and (\ref{ab10}) give pairwise distinct $\widehat d_j \in \C$ with 
$$
\frac{f_j'(z)}{f_j(z)} + \frac{v_2(z)}{u_2(z)} = \widehat d_j z^{\rho-1} + \ldots .
$$
If $\widehat d_j \neq 0$ then $f_j' u_2 + f_j v_2$  is again a canonical formal solution of $M[y]=0$ with exponential part~$q_j$,
and  so a constant multiple of $w_j$.
Since $k \geq 4$,  this implies in view of (\ref{ab3}) that
\begin{equation}
 \label{ab14}
\mu_j = \lambda_j + \eta_2 + \rho - 1 
\end{equation}
for $j=1$ and $j=2$, or for $j=k-1$ and $j=k$.
If (\ref{ab14}) holds for $j=1$ and $j=2$ then (\ref{ab8}), (\ref{ab9}) and (\ref{ab14}) give 
$$
\mu_1 = \lambda_1 + \eta_2 + \rho - 1 = \lambda_2 + \eta_1 + \rho- 1 ,
\quad
\mu_2 =  \lambda_2 + \eta_2 + \rho - 1 = \lambda_1 + \eta_3 + \rho- 1 ,
$$
from which it follows that 
$$
\eta_1 - \eta_2 = \lambda_1 - \lambda_2 = \eta_2 - \eta_3 , \quad \eta_1 + \eta_3 - 2 \eta_2 = 0.
$$
Similarly, if (\ref{ab14}) holds for $j=k-1$ and $j=k$,
then (\ref{ab8}), (\ref{ab9}) and (\ref{ab14}) give 
$$
\mu_k = \lambda_k + \eta_2 + \rho - 1 = \lambda_{k-1} + \eta_3 + \rho- 1 ,
\quad
\mu_{k-1} =  \lambda_{k-1} + \eta_2 + \rho - 1 = \lambda_k + \eta_1 + \rho- 1 ,
$$
which delivers 
$$
\eta_1 - \eta_2 = \lambda_{k-1} - \lambda_k = \eta_2 - \eta_3 , \quad  \eta_1 + \eta_3 - 2 \eta_2 = 0.
$$
\hfill$\Box$
\vspace{.1in}

\begin{lem}
 If $b_3 = 0$ in (\ref{ac1}) then $f$ satisfies the conclusion of Proposition \ref{mainprop}.
\end{lem}
\textit{Proof.} 
Using (\ref{5}) write, on $\Sigma$,
\begin{equation}
 \label{D3a1}
g = U_1 + U_2, \quad - \frac{f'}{f} = \frac{V[g]}g =
\frac{V[U_1] + V[U_2]}{U_1+U_2} = \frac{ V[U_1]/U_2 + V[U_2]/U_2}{ e^\Phi + 1}, \quad e^\Phi = \frac{U_1}{U_2}  .
\end{equation}
A zero of $g$ arises wherever $U_1/U_2 = e^\Phi = -1$, and the multiplicity of the pole of $f$ at such a point is 
\begin{equation}
 \label{D3a2}
m_0 = \frac{ V[U_1]/U_1 - V[U_2]/U_2}{  \Phi' } .
\end{equation}
By (\ref{ab2}) and (\ref{ac1}), the function 
$\zeta = (1/\pi i) \Phi = (1/ \pi i) \log U_1/U_2$ maps the sector 
$\Sigma$ univalently onto a region containing a half-plane $\pm {\rm Re} \, \zeta > M_1 \in \R$,  and 
(\ref{D3a2}) holds wherever $\zeta$ is an odd integer. Thus (\ref{ab2}), (\ref{ac1}), (\ref{ab7}) and Lemma \ref{lemintegervalued} give a polynomial $Q^*$ such that 
\begin{equation}
 \label{D3a3}
\frac{V[U_1]}{U_1} - \frac{V[U_2]}{U_2} =  Q^*(\Phi) \Phi' , \quad U_2 V[U_1] - U_1 V[U_2] =  Q^*(\Phi) ( U_2U_1' - U_1 U_2'). 
\end{equation}

Suppose first that $Q^*(\Phi)$ is rational at infinity in (\ref{D3a3}). Then it follows from Lemma \ref{newdelem}
that $U_1$ and $U_2$ solve a second order equation (\ref{7}) with $E_1$ and $E_0$ rational at
infinity, and so does $g$, by (\ref{D3a1}), contradicting Lemma \ref{lemmainstep}(C).

It may therefore be assumed henceforth that  $Q^*$ is non-constant. 
Then (\ref{D3a2}) and (\ref{D3a3}) show that the multiplicity $m_0(z)$ of a pole $z \in \Sigma$ of 
$f$  tends to $\infty$ as $z \to \infty$, faster than  $|z|^{\rho_1}$ for some $\rho_1 > 0$. Since the zeros of $g = U_1+U_2$ in $\Sigma$ have of exponent of convergence
$\rho$, this is incompatible with Case B of Lemma \ref{sofar}. Hence 
Case A of Lemma \ref{sofar} must hold, and so $(\eta_1 - \eta_2 ) - ( \eta_2 - \eta_3 ) = \eta_1 + \eta_3 - 2 \eta_2 $ is a positive integer. 

Furthermore, the left-hand side of (\ref{D3a3}) has a meromorphic 
continuation along any path in $\Omega(r_1)$, as has $\Phi'$,
but if a continuation of $U_1/U_2$ has a zero or pole at some $z_0$ then $\Phi(z) = \log U_1(z)/U_2(z) $ behaves like
$m_1 \log (z-z_0) $ as $z \to z_0$, for some $m_1 \in \Z \setminus \{ 0 \}$. Therefore  (\ref{D3a3}) implies that $e^\Phi = U_1/U_2$ continues without poles
or zeros in $\Omega(r_1)$, and so any zeros of continuations of $U_1$ and $U_2$ are shared.  

Take any sector $\Sigma^*$ given by $| \arg z - \theta^* | \leq \pi /  \rho - \delta_1 $,  where 
${\rm Re} \, P(re^{i\theta^*} ) = O( r^{\rho-1})$ as $r \to \infty$,
let $\widetilde U_1, \widetilde U_2 $ be  continuations of $U_1, U_2$ to $\Sigma^*$, and write 
\begin{equation}
 \label{newUj}
\widetilde U_1 = d_1 \psi_1 + d_2 \psi_2 + d_3 \psi_3 , \quad 
\widetilde U_2 = e_1 \psi_1 + e_2 \psi_2 + e_3 \psi_3 , \quad d_j, e_j \in \C, 
\end{equation}
on $\Sigma^*$, in which the $\psi_j$ are analytic solutions of $N[y]=0$ which satisfy, as $z \to \infty$ on $\Sigma^*$,
\begin{equation}
 \label{newUj2}
\psi_1 (z) = z^{\eta_1} e^{P(z)} (1+o(1)), \quad \psi_2 (z) = z^{\eta_2} (1+o(1)), \quad \psi_3 (z) = z^{\eta_3} e^{-P(z)} (1+o(1)). 
\end{equation}
Suppose that $\widetilde U_1$ and $\widetilde U_2$ have a sequence $\zeta_\mu \to \infty$ of common zeros in $\Sigma^*$. 
The matrix with rows $(d_1, d_2, d_3)$ and $(e_1, e_2, e_3)$ has rank $2$, since $U_1$ and $U_2$ are linearly independent, and so Cramer's rule gives
$e_4, e_5 \in \C$ and a permutation $(j, j', j'')$ of $(1, 2, 3)$  such that 
$$
\psi_{j'} (\zeta_\mu) = e_4 \psi_j(\zeta_\mu), \quad 
\psi_{j''} (\zeta_\mu) = e_5 \psi_j(\zeta_\mu) \quad \hbox{as $\mu \to \infty$.}
$$
Here $e_4e_5 \neq 0$, as $\psi_j (\zeta_\mu) \neq 0$ for large $\mu$. But this gives a contradiction, 
since the fact that $(\eta_1 - \eta_2 ) - ( \eta_2 - \eta_3 )$ is positive implies that 
$\psi_2(\zeta_\mu)/\psi_3(\zeta_\mu)= o( |\psi_1(\zeta_\mu)/\psi_2(\zeta_\mu)|)$ as $\mu \to \infty$.

It follows that $U_1$ and $U_2$ continue without zeros in some annulus $\Omega(r^*)$. 
Lemma \ref{lemPL} shows that there exists $\rho_2 > 0$ such that any continuation of $U_2$ to any sector in $\Omega(r^*)$ satisfies 
$\log |U_2(z)| = O\left( |z|^{\rho_2} \right)$
as $z \to \infty$ there.
Take a sector $\Sigma^{**}$ given by $\theta_1 < \arg z < \theta_2$, where these $\theta_j$ are such that no $\theta \in [ \theta_1, \theta_2 ]$ has 
$\mathrm{Re} \, P(re^{i \theta}) = O( r^{\rho-1} )$ as $r \to \infty$. 
For any continuation of $U_2$ to $\Sigma^{**}$ there exist
$P^* \in \{  -P, 0, P \}$ and a matching $\eta^* \in \{ \eta_1, \eta_2, \eta_3 \}$ such that $U_2(z) \sim c z^{\eta^*} \exp( P^*(z))$ as $z \to \infty$ in $\Sigma^{**}$. 
Since $U_2(z) \sim c z^{\eta_2}  $ as $z \to \infty$ in $\Sigma$, repeated application of the 
Phragm\'en-Lindel\"of principle to the continuations of $U_2(z) z^{-\eta_2}$ or its reciprocal shows that $P^* =0$,
and so $\eta^* = \eta_2$.  Examining (\ref{newUj}) in the light of (\ref{newUj2}), first on a subsector of $\Sigma^*$ on which $e^P$ is large
and subsequently on a subsector where $e^{-P}$ is large, 
forces $0 = e_1 = e_3$.  
Choosing $\Sigma^* = \Sigma$ gives $e_0 \in \C$ such that $z^{e_0}U_2(z)$ is analytic and zero-free of finite order of growth in some annulus $\Omega(r^{**})$.
This, coupled with almost identical reasoning applied to $U_1$, 
shows that $U_1'/U_1$,  $U_2'/U_2$ and $\Phi'$ are rational at infinity, as is $Q^*(\Phi)$ by (\ref{D3a3}), and this case has already been dealt with. 
\hfill$\Box$
\vspace{.1in}

Assume henceforth that $b_1b_3 \neq 0$ in (\ref{ac1}), and write this formula for $g$ as 
\begin{equation}
\label{ac2}
 g = A e^{-P} ( (e^P - B)^2 - C^2) ,
\quad U_1 = Ae^P, \quad U_2 = -2 AB, \quad U_3 = A( B^2 - C^2 )e^{-P} . 
\end{equation}
By (\ref{ab2}), this initially formal expression for $g$ results in, as $z \to \infty$ in $\Sigma$, 
\begin{eqnarray}
  \label{ac3}
A(z) &=&  b_1 z^{\eta_1} \chi_1 (z) , \quad 
B(z) = - \frac{b_2}{2b_1} z^{\eta_2 - \eta_1} \chi_2 (z) , \nonumber \\
\quad B(z)^2-C(z)^2 &=&  \frac{b_3}{b_1} z^{\eta_3 - \eta_1} \chi_3(z), 
\quad \chi_j(z) = 1 + o(1)  .
\end{eqnarray}
Here the $\chi_j$  have asymptotic series on $\Sigma$ in descending integer powers of $z$ and,  by Lemma \ref{sofar},
$\eta_3 - \eta_1 - 2( \eta_2 - \eta_1)  =  \eta_1 + \eta_3 - 2 \eta_2  $ is a non-negative even integer. 
Evidently $A, B$ and $ E=C^2$ are analytic on $\Sigma$, and $E$ does not vanish identically, 
since zeros of $g$ are simple.
Furthermore, it is clear from (\ref{ac2}) 
that, at a zero of $g$ in $\Sigma$, 
\begin{equation}
 \label{ac6}
(e^P - B)^2 =   E = C^2, \quad g' = Ae^{-P} ( 2(e^P - B) (P'e^P - B') - E' ) .  
\end{equation}

\begin{lem}
\label{lemCan}
Let $d = \pm 1$. Then there exist $r_2 > 0$ and 
$\sigma_d , \tau_d \in \C^* $, as well as  $ \gamma_d, \zeta_d \in \C $, such that $B+dC$ is analytic on $\Sigma \cap \Omega(r_2)$ and
\begin{equation}
 \label{ac4}
C(z)^2 = \sigma_d  z^{\gamma_d} \psi_1(z) , \quad \psi_1(z) = 1+o(1) , 
\end{equation}
and 
\begin{equation}
 \label{ac8}
B(z) +d C(z) =  \tau_d z^{\zeta_d} \psi_2(z), \quad \psi_2(z) = 1+o(1) ,
\end{equation}
as $z \to \infty$ in $\Sigma$, 
in which  the $\psi_j(z) $  have asymptotic series on $\Sigma$
in descending integer powers of $z^{1/2}$. Furthermore, if conclusion (A) of Lemma \ref{sofar} holds, 
then 
$\gamma_d = \eta_3 - \eta_1$. 
 \end{lem}
\textit{Proof.} Note first that $B+dC$ does not vanish identically, since $B^2-C^2$ does not. 
All conclusions of the lemma clearly follow
from (\ref{ac3}) if $b_2 = 0$ or $\eta_3 - \eta_1 - 2( \eta_2 - \eta_1) > 0$, and in particular if conclusion (A) of Lemma \ref{sofar} holds.

Assume therefore that $b_2 \neq 0$ and $\eta_3 - \eta_1 = 2( \eta_2 - \eta_1)$. 
Then (\ref{ac3}) implies that $\widetilde C (z) = C(z)^2 z^{2(\eta_1-\eta_2)} $ has an asymptotic series on $\Sigma$ in descending non-positive integer powers of $z$. 
If this asymptotic series for $\widetilde C(z)$  vanishes identically then, by making $\Sigma$ slightly narrower
if necessary, it may be assumed that 
$E(z) = C(z)^2$ and $E'(z)$ both tend
to zero in $\Sigma$ transcendentally fast, that is,
faster than any negative power of $z$, but $f$ still has infinitely many poles there.
This implies using (\ref{ac7}) and (\ref{ac6})
that if $M_1$ is a positive integer and $z$ is a pole of $f$ of multiplicity $m_0(z)$ in 
$\Sigma$,  with $|z|$ large, then
$$
g(z) = 0, \quad e^{P(z)} =  B(z) + O\left( |z|^{-2M_1} \right), \quad g'(z) =  O\left( |z|^{-M_1} \right), \quad |z|^{M_1} = O( m_0(z)) ,
$$
which is a contradiction since $f$ has finite order. 
Hence there must exist an integer $m_1 \leq 0$ such that (\ref{ac4}) holds with $\gamma_d = 2( \eta_2 - \eta_1) + m_1$,
in which   $\psi_1(z)$  has  an asymptotic series in descending integer powers of $z$.
It is now clear from (\ref{ac3}) and (\ref{ac4}) that $\widetilde B (z) = (B(z)+dC(z))z^{\eta_1-\eta_2} $ has an asymptotic series 
on $\Sigma$ in descending  integer powers of $z^{1/2}$; thus (\ref{ac8}) holds unless 
this series for $\widetilde B (z)$ vanishes identically, in which case $B(z)+dC(z)$  tends to zero transcendentally fast on $\Sigma$, and so does 
$B(z)^2 - C(z)^2$, by the second equation of (\ref{ac3}),  which forces $b_3=0$ in (\ref{ac3}), contrary to assumption.
\hfill$\Box$
\vspace{.1in}

\begin{lem}
 \label{lemQd}
For $d = \pm 1$ 
there exists a polynomial $Q_d \not \equiv 0$ such that 
\begin{equation}
 \label{ac9}
\left[ 2d CA \left( P' - \frac{B'+dC'}{B+dC} \right) \right]^{-k} = Q_d (P - \log (B+dC) ) .
\end{equation} 
\end{lem}
\textit{Proof.} 
The
function $g$ has a zero in $\Sigma$  wherever $e^P = B + dC$,  and at such a zero (\ref{ac6}) gives
\begin{eqnarray}
 g' &=& Ae^{-P} ( 2dC (P'e^P - B') - 2CC' ) = 2dCA e^{-P} ( P'e^P - B' -dC')  \nonumber \\
&=&  2dCA \left( P' - \frac{B'+dC'}{B+dC} \right) .
\label{ac9aaa}
\end{eqnarray}
Here  (\ref{ac8})  shows that 
$\zeta = (1/2\pi i) (P(z) - \log (B(z)+dC(z)))$ maps 
a subdomain of $\Sigma$ univalently onto 
a half-plane $\pm {\rm Re} \, \zeta > M_1 \in \R$. 
Because (\ref{ac7}) implies that $(g')^{-k}$ is 
integer-valued at each zero of $g$, and so at points where $\zeta$ is integer-valued, it follows from (\ref{ac3}), (\ref{ac4}) and
Lemma \ref{lemintegervalued} that 
a polynomial $Q_d$ exists as asserted.
\hfill$\Box$
\vspace{.1in}

\begin{lem}
 \label{lemQd2a}
For  $d = \pm 1$ the polynomial 
$Q_d$ in (\ref{ac9}) is constant. 
\end{lem}
\textit{Proof.}
Assume that
$Q_d$ is non-constant. 
Then it follows from (\ref{ac7}), (\ref{ac8}),  (\ref{ac9})  and (\ref{ac9aaa})
that the multiplicity $m_0(z)$ of the pole of $f$ at $z \in \Sigma$ tends to $\infty$ faster than some positive power of $|z|$ 
and, since the exponent of convergence of the zeros of $e^P - (B+dC)$ in $\Sigma$ is $\rho$, 
this implies that $N(r, f)$ has order greater than $\rho$, which is incompatible with 
conclusion (B) of Lemma \ref{sofar}. 

Hence conclusion (A) of Lemma \ref{sofar} must hold. In view of  (\ref{ac3}) and Lemma~\ref{lemCan},  it follows that 
$\eta_1 + \eta_3 = - 2 (\rho-1)$ 
and $\gamma_d  = \eta_3 - \eta_1$, and that 
$$C(z)A(z) \sim c z^{\gamma_d/2 + \eta_1 } = c z^{(\eta_1+ \eta_3)/2 } = c z^{ 1 - \rho } $$
as $z \to \infty$ in $\Sigma$. But then the left-hand side of (\ref{ac9}) is bounded as $z \to \infty$ in $\Sigma$,
which is a contradiction.  
\hfill$\Box$
\vspace{.1in}

\begin{lem}
\label{lemBCconst}
There exist a large positive $r_3 $ and an analytic function $K$ such that 
\begin{equation}
\label{U2U3}
 K' = \frac1{U_1}, \quad U_1 = Ae^P, \quad U_2 = -2 AB = e_3 U_1 K , \quad U_3 = A(B^2-C^2)e^{-P} = e_4 U_1 K^2 , 
\end{equation}
on $\Sigma \cap \Omega(r_3) $, where $e_3, e_4 \in \C$ and $e_4 \neq 0$. 
\end{lem}
\textit{Proof.} Suppose first that $B \not \equiv 0$. Then (\ref{ac9}) holds for $d=1$ and $d=-1$, with $Q_1$ and $Q_{-1}$ both constant
by Lemma \ref{lemQd2a}. Hence, by (\ref{ac8}) and (\ref{ac9}), 
$$
 CA \left( P' - \frac{B'+C'}{B+C} \right), \quad CA \left( P' - \frac{B'-C'}{B-C} \right)
$$
are  both  constant, and so identically equal.
Thus  $(B+C)/(B-C)$ must be constant and so must $B/C$.
Now (\ref{ac2}), (\ref{ac4}), (\ref{ac9}) and Lemma \ref{lemQd2a}
yield, with  $c \in \C^*$ as before, 
\begin{equation}
 \label{U2U3b}
CA \left( P' - \frac{C'}{C} \right) = c, 
\quad C' - P'C = \frac{c}A, \quad C e^{-P} = c \int \frac{1}{Ae^P} = c \int \frac{1}{U_1} = cK ,
\end{equation}
from which (\ref{U2U3}) follows, using (\ref{ac2}) again. On the other hand, if $B \equiv 0$ then the first equation of 
(\ref{U2U3b}) still holds, by Lemma \ref{lemQd2a}, and the formula for $U_2$ in (\ref{U2U3}) is trivially satisfied with $e_3=0$. 
\hfill$\Box$
\vspace{.1in}

\begin{lem}
\label{lemac}
The function $K$ of Lemma \ref{lemBCconst}
continues meromorphically along any path in the annulus $\Omega(r_3)$, its continuations locally univalent. 
Moreover, all  zeros of any continuation of $U_1$ into $\Omega(r_3)$ are 
simple poles of $K$.
\end{lem}
\textit{Proof.} Since $e_4 \neq 0$ in (\ref{U2U3}), writing
\begin{equation}
 \label{Heqn}
\Phi = K^2 = \frac{U_3}{e_4U_1},  \quad 
 \frac1{U_1^2} = (K')^2 = \frac{(\Phi')^2}{4\Phi} ,
\end{equation}
shows that $\Phi$  continues meromorphically along any path in  $\Omega(r_3)$.
Any zero of any
continuation of $U_1$ is either simple or double, since $U_1$ solves $N[y]=0$, and must be a pole of $\Phi$, by (\ref{Heqn}). 
Comparing multiplicities in (\ref{Heqn}) excludes simple zeros of $U_1$, and double zeros of $U_1$ have to be triple poles of $\Phi'$ 
and so double poles of~$\Phi$. Furthermore, any zeros of any continuation of $\Phi$ must be double, again by (\ref{Heqn}). Thus $K = \Phi^{1/2}$ continues meromorphically
along paths in $\Omega(r_3)$, and is locally univalent since $K' (z) = 1/U_1 (z) \neq 0$. 
\hfill$\Box$
\vspace{.1in}

Again because $e_4 \neq 0$ in (\ref{U2U3}), there exists a polynomial $Q_2$ of degree $2$ such that (\ref{ac1}) and continuation of $g$
into $\Omega(r_3)$ give
$g = 
Q_2(K)/K'$, whether or not $U_2 \equiv 0$,
where $K$ is as in Lemma \ref{lemac}. Hence $g = 0$ forces $K = a$, where $Q_2(a) = 0$, and so $g'= Q_2'(a) = \pm b$ for some $b \in \C^*$, by elementary
properties of
quadratics.
It now follows using (\ref{ac7}) 
that all  poles of $f$ in $\Omega(r_3)$ have the same multiplicity,
and $f$ satisfies the conclusion of Proposition \ref{mainprop}. 
\hfill$\Box$
\vspace{.1in}

\section{The case of a repeated trivial exponential part}\label{repeattriv}

There remains only one case to deal with, in which
the equation $N[y]=0$ has two linearly independent formal solutions 
$g_1$, $g_2$ with trivial exponential part.
The third exponential part $\kappa$ must be non-zero, by Lemma \ref{lemkappanotsame}.
Take 
a ray $\arg z = \theta_0$ 
on which
$\kappa \prec 0$, and 
label the exponential parts arising from $L$ and $M$ to be consistent with (\ref{yz2})
on $\arg z = \theta_0$. 
It then follows from  Lemma \ref{prop2a} that 
the exponential parts $q_j $ for the equation
$L[y]=0$ are pairwise distinct, and the same is true for $M[y]=0$, and so the 
formal solutions
of these equations are free of logarithms. This implies that 
any formal solution $G$ of $N[y]=0$ 
is also free of logarithms; to see this, take a fundamental set of canonical formal solutions $f_j$ 
of $L[y]=0$, write $f_j'G + f_j V[G] = w_j$, 
where the $w_j$ are formal solutions of $M[y]=0$, and solve for $G$.

Therefore $N[y]=0$ has  linearly independent canonical 
formal solutions $g_1$, $g_2$ each having the form $g_j(z) = z^{m_j} ( 1 + \ldots )$, with $m_j \in \C$.
There exists a  third  canonical formal solution $g_3$, which has  exponential part $\kappa$ and,
by Lemma \ref{lemexppart}, annihilates some canonical formal solution $h_\mu$ of $L[y]=0$, with exponential part $q_\mu$ say. 
Consider the terms $R_j = V[g_j]/g_j$, for $j=1, 2$; these are 
formal series  in descending integer powers of $z$.
Hence   
$$ S_j =  h_\mu' g_j + h_\mu V[g_j] = h_\mu g_j \left(h_\mu'/h_\mu +  R_j   \right) $$
is a  formal
solution of $M[y] =0$, for $j=1, 2$, and either is identically zero or has exponential 
part $q_\mu$. Since the exponential parts for $M$ are all different,
$S_1$ and $S_2$ must be linearly 
dependent, and some non-trivial linear combination $g_4$ of $g_1$ and $g_2$ must annihilate
$h_\mu$, as does $g_3$. 
Therefore
$g_3 V[g_4] = g_4 V[g_3]$ and 
Lemma
\ref{d=0lemA}, with $d=0$,  gives an equation (\ref{7}) solved by $g_3$ and $g_4$.
Furthermore, $g_4$ must be a canonical formal solution of $N[y]=0$; this is obvious unless $g_4 = d_1 g_1 - d_2 g_2$ with
$d_1, d_2 \in \C^*$, in which case 
$$d_1 S_1 = d_2 S_2, \quad
d_1 g_1 \left(h_\mu'/h_\mu +  R_1   \right) = d_2 g_2 \left(h_\mu'/h_\mu +  R_2   \right).
$$
Thus $W(g_3, g_4)$ has non-zero exponential part, 
so that $E_1(\infty) \neq 0$ in (\ref{7}), and the conclusion of Proposition~\ref{mainprop} follows from Lemma \ref{d=0lemA}.
\hfill$\Box$

\noindent
School of Mathematical Sciences, University of Nottingham, NG7 2RD, UK.\\
james.langley@nottingham.ac.uk

\end{document}